\documentclass{article}
\usepackage{float}
\usepackage{amsmath}
\newtheorem{theorem}{Theorem}
\newtheorem{lemma}{Lemma}

\usepackage{amssymb}
\usepackage{graphicx,bm}
\usepackage{url}
\usepackage[colorlinks,linkcolor=blue,citecolor=green,urlcolor=green]{hyperref}
\usepackage{color}
 \usepackage[compress]{cite}
 \usepackage{enumerate}
 \usepackage{epstopdf}
 \usepackage{multirow}
 
  \numberwithin{figure}{section}
    \numberwithin{theorem}{section}
    \numberwithin{remark}{section}
  \numberwithin{table}{section}
  \newcommand{\Ii}{{\mathrm{i}}}
  \numberwithin{equation}{section}

\newcommand{\red}[1]{{\color{red}#1}}

\usepackage[margin=1in]{geometry}
 \newcommand{\T}{\mathsf{T}}

 \newcommand{\CP}{{{\bm P}}}

 \newcommand{\IC}{{\mathbb{C}}}

 \newcommand{\tn}[1]{\ \textnormal{#1}\ }
 \newcommand{\bmt}{\left[ \begin{array}{ccccccccccccccccccccccccccccccccccccc}}
 	\newcommand{\emt}{\end{array}\right]}
 \usepackage{booktabs}

\makeatletter
\renewcommand\normalsize{%
   \@setfontsize\normalsize\@xpt\@xiipt
   \abovedisplayskip 6.8\p@ \@plus2\p@ \@minus5\p@
   \abovedisplayshortskip \z@ \@plus3\p@
   \belowdisplayshortskip 6.8\p@ \@plus3\p@ \@minus3\p@
   \belowdisplayskip \abovedisplayskip
   \let\@listi\@listI}
\makeatother

\begin{document}

 \title{{\bf ParaDiag}: Parallel-in-Time Algorithms \\Based on {the}
   {\em Diagonalization} Technique}

\author{ Martin J. Gander$^{\color{red}{a}}$ ~~~   Jun Liu$^{\color{red}{b}}$   ~~~ Shu-Lin Wu$^{\color{red}{c}}$  ~~~  Xiaoqiang Yue$^{\color{red}{d}}$~~~ Tao Zhou$^{\color{red}{e}}$\\
$^{\color{red}{a}}${\em  Section of Mathematics, University of Geneva, CH-1211 Geneva, Switzerland}\\
E-mail: \texttt{martin.gander@unige.ch}\\
$^{\color{red}{b}}${\em Department of Mathematics and Statistics},\\ {\em Southern Illinois University Edwardsville},
{\em Edwardsville, IL 62026, USA} \\
E-mail:  \texttt{juliu@siue.edu}\\
$^{\color{red}{c}}${\em School of Mathematics and Statistics, Northeast Normal University},\\
{\em Changchun 130024, China}\\
 E-mail: \texttt{wushulin84@hotmail.com}\\
 $^{\color{red}{d}}${\em School of Mathematics and Computational Science, Xiangtan University},\\
 {\em Xiangtan 411105, China}\\
 E-mail: \texttt{yuexq@xtu.edu.cn}\\
 $^{\color{red}{e}}${\em Institute of Computational Mathematics and Scientific/Engineering Computing},\\
 {\em AMSS, the Chinese Academy of Sciences, Beijing, China}\\
  Email: \texttt{tzhou@lsec.cc.ac.cn}}

\maketitle

\tableofcontents

\begin{abstract}
In 2008, Maday and R{\o}nquist introduce{d} an interesting {new
  approach for the} direct parallel-in-time (PinT) {solution} of
time-dependent PDEs. The idea {is to diagonalize the time stepping
  matrix}, keeping the matrices for the space discretization
unchanged{, and then to solve all time steps in parallel}. Since
{then, several variants appeared}, and we {call these closely
  related algorithms} {\em ParaDiag} {algorithms}. ParaDiag
algorithms {in the literature} can be {classified into  {two}
  groups}:
\begin{itemize}
\item ParaDiag-I: direct {standalone} solvers,
\item ParaDiag-II:  iterative {solvers}{.}
 \end{itemize}
{We will} explain the basic feature{s} of each {group} in
this note. {To have concrete examples, we} will introduce
ParaDiag-I and ParaDiag-II for the advection-diffusion equation. {We will also introduce   ParaDiag-II  for the wave equation  and an optimal control
{problem for the} wave equation}. {We could have used the
  advection-diffusion equation as well to illustrate} {ParaDiag-II},
{but wave equations are known to cause problems for certain PinT
  algorithms and thus constitute an especially interesting example for
  which ParaDiag algorithms were tested}. We show the main {known}
theoretical results {in each case, and also provide Matlab codes}
for testing.  The goal of the Matlab codes is to help the interested
reader understand the key {features} of the {ParaDiag}
algorithms, {without intention to be highly tuned for efficiency and/or
  low} memory {use}.

We also provide {speedup measurements} of ParaDiag {algorithms} for a
2D linear advection-diffusion equation.  These results are obtained on
the Tianhe-1 supercomputer in China {and the SIUE Campus Cluster
  in {the} US}, which is a multi-array, configurable and cooperative
parallel system{, and we compare these results to the performance of
  parareal and MGRiT, two widely used PinT algorithms}.  In {a}
forthcoming update of this note, we will {provide} more material {on
  ParaDiag algorithms, in particular further} Matlab codes and
parallel comput{ing} results, {also for more realistic applications}.
\end{abstract}

\section{Basic idea of ParaDiag}

{We start with a basic introduction to ParaDiag algorithms.}
Suppose we need to {solve in parallel the system of ODEs}
$M\dot{U}(t)+KU(t)=f(t)$ with initial value $U(0)=U_0$ arising from
the semi-discretization of a time-dependent PDE, where
$M,K\in\mathbb{C}^{N_x\times N_x}$.  For finite element
discretizations, $M$ {is the mass matrix} and $K$ is the stiffness
matrix. For finite difference discretizations, $M=I_x$ {is just an
  identity matrix. The classical approach for} solving {such
  systems of ODEs} is to apply {a} time-integrator, and then solve
the resulting difference equation step-by-step {in time}.  Instead,
ParaDiag tries to solve these difference equations {\em all-at-once}.
For linear multi-step methods, the all-at-once system is {of the
  form}
\begin{equation}\label{eq1.1}
{\bm A}{\bm u}={\bm b}, ~{\bm A}:=B_1\otimes M+B_2\otimes K,
\end{equation}
where $B_1, B_2\in\mathbb{R}^{N_t\times N_t}$ are Toeplitz matrices
specified by the time-integrator and $N_t$ is the number of time
steps\footnote{For Runge-Kutta methods, the all-at-once system is
  different and will be treated in a forthcoming update of this
  note.}.  All ParaDiag algorithms {focus on} treating the matrices
$B_1$ and $B_2$, while keeping $M$ and $K$ unchanged.  There are
mainly two approaches. First, by rewriting {problem}
\eqref{eq1.1} as
\begin{equation}\label{eq1.1a}
{\bm A}{\bm u}=(B_2^{-1}\otimes I_x){\bm b}, ~{\bm A}:=B_2^{-1}B_1\otimes M+I_t\otimes K,
\end{equation}
we directly diagonalize the matrix $B=B_2^{-1}B_1$.  This leads to
ParaDiag-I, the group of direct PinT solvers.  {The research on
  ParaDiag-I focuses on {obtaining time stepping matrices} $B$
  {that are} diagonalizable {with a} condition number of the
  {associated} eigenvector matrix as small as possible. {A
    c}oncrete example is the {original} algorithm} based on using
different {time} step sizes $\{\Delta t_n\}$, e.g., a geometrically
increasing sequence $\Delta t_n=\Delta t_1\tau^{n-1}$ with $\tau>1$,
{which makes} the time-discretization matrices diagonalizable
\cite{MR08,GH16,GH19}.  {{We will show} new progress on
    ParaDiag-I in Section \ref{sec2.1.2}, {using} a hybrid
    time-discretization {with uniform time} step sizes.  }

The second treatment is to use a uniform step size $\Delta t$ and
solve the all-at-once system \eqref{eq1.1} iteratively{, which
  leads to ParaDiag algorithms in the ParaDiag-II
  group.} There are {several variants}, but the {common point is
  to introduce} the $\alpha$-circulant block matrix
\begin{equation}\label{eq1.2}
 {\bm P}_\alpha:=C_1^{(\alpha)}\otimes M+C_2^{(\alpha)}\otimes K,
\end{equation}
where $C_1^{(\alpha)}$ and $C_2^{(\alpha)}$ are Strang type
$\alpha$-circulant matrices {constructed from} $B_1$ and $B_2$, and
$\alpha\in(0, 1]$ is a free parameter.  One can then either solve
  \eqref{eq1.1}
  via the stationary iteration \cite{WL20b}
\begin{equation}\label{stationaryiteration}
{\bm P}_\alpha {\bm u}^k=({\bm P}_\alpha-{\bm A}){\bm u}^{k-1}+{\bm b},
\end{equation}
where $k\geq1$ is the iteration index, or via Krylov subspace methods (e.g., GMRES, MINRES) by
solving the preconditioned system \cite{MPW18}
\begin{equation}
  {\bm P}_\alpha^{-1}{\bm A} {\bm u}={\bm P}_\alpha^{-1}{\bm b},
\end{equation}
which is nothing else than the stationary iteration
\eqref{stationaryiteration} written at its fixed point, i.e. at
convergence.

The algorithms proposed in \cite{W18} and \cite{GW19} are essentially
 ParaDiag-II algorithms {as well}, but they are derived from a different
point of view.  For example, in \cite{GW19} the authors introduced a
{Waveform Relaxation} (WR) iteration $M\dot{U}^k(t)+KU^k(t)=f(t)$,
$U^k(0)=\alpha (U^k(T)-U^{k-1}(T))+U_0$, and after a
time-discretization {one} can show that at each iteration the
all-at-once system is ${\bm P}_\alpha {\bm u}^k={\bm b}^{k-1}$, where
${\bm b}^{k-1}=({\bm P}_\alpha-{\bm A}){\bm u}^{k-1}+{\bm b}$.  The
algorithm in \cite{W18} can be understood similarly.

For {each variant  of   ParaDiag-II  we need to compute} ${\bm
  P}^{-1}_\alpha{\bm r}$ with ${\bm r}$ being an input vector.  The
reason for using ${\bm P}_\alpha$ is twofold{: first}, since
$C_1^{(\alpha)}$ and $C_2^{(\alpha)}$ are Strang type
$\alpha$-circulant matrices {constructed from the Toeplitz
  matrices} $B_1$ and $B_2$, it naturally holds that ${\bm P}_\alpha$
converges to ${\bm A}$ as $\alpha$ goes to zero. This implies that by
using a relatively small $\alpha$, the ParaDiag-II
algorithms converge rapidly. The second point lies in the fact that
$C_1^{(\alpha)}$ and $C_2^{(\alpha)}$ can be diagonalized
simultaneously, as {is shown in the following Lemma.}
 \begin{lemma}[see \cite{B05}]\label{lem1}
Let
$\mathbb{F}=\frac{1}{\sqrt{N_t}}\left[\omega^{(l_1-1)(l_2-1)}\right]_{l_1,
  l_2=1}^{N_t}$ (with $\Ii=\sqrt{-1}$ and
$\omega=e^{\frac{2\pi{\Ii}}{N_t}}$) be the discrete Fourier matrix and
define {for any given parameter $\alpha\in (0,1]$ the} diagonal
  matrix
\begin{equation*}
\Gamma_\alpha=\begin{bmatrix}1 & & &\\ &\alpha^{\frac{1}{N_t}} & &\\
& &\ddots &\\ & & &\alpha^{\frac{N_t-1}{N_t}}\end{bmatrix}.
\end{equation*}
{Then the} two $\alpha$-circulant matrices $C^{(\alpha)}_1,
C^{(\alpha)}_2\in\mathbb{C}^{N_t\times N_t}$ can be simultaneously
diagonalized as
$$
C^{(\alpha)}_{j}=V D_{j}V^{-1},~ D_{j}={\rm diag}\left(\sqrt{N_t}\mathbb{F}\Gamma_\alpha C^{(\alpha)}_{j}(:,1)\right), ~j=1, 2,
$$
where $V=\Gamma_\alpha^{-1}\mathbb{F}^*$ and $C^{(\alpha)}_{j}(:,1)$
{represents} the first column of $C^{(\alpha)}_{j}$, $j=1,2$.
\end{lemma}
Due to the property of {the} Kronecker product, we can
{factor} ${\bm P}_\alpha=(V\otimes I_x)(M\otimes
D_1+A\otimes D_2)(V^{-1}\otimes I_x)$ and thus we can compute ${\bm
  P}^{-1}_\alpha{\bm r}$ by {performing the following} three steps:
\begin{equation}\label{eq1.3}
\begin{split}
&\text{Step-(a)} ~~S_1=(V^{-1}\otimes I_x){\bm r},\\
&\text{Step-(b)} ~~S_{2, n}=\left(\lambda_{1,n}M+\lambda_{2,n}A\right)^{-1}S_{1,n}, ~n=1,2,\dots,N_t,\\
&\text{Step-(c)} ~~{\bm u}=(V\otimes I_x)S_2,
\end{split}
\end{equation}
where $S_1=(S_{1,1}^\top,\dots, S_{1,N_t}^\top)^\top$ and
$S_2=(S_{2,1}^\top,\dots, S_{2,N_t}^\top)^\top$.  Since $V$ and
$V^{-1}$ are given {by FFT techniques, Step-(a) and Step-(c) can
  be} computed efficiently with $O(N_xN_t\log N_t)$
operations. Step-(b) can be computed in parallel since {all linear
  systems} are completely independent from each other at different
time points.  These three steps {represent the key steps of
  ParaDiag algorithms and will appear frequently in this note, although the
  details differ in the various cases}.

For nonlinear problems $M\dot{U}+f(U)=0$ with $U(0)=U_0$, the basic
idea for {applying ParaDiag algorithms} is as follows: for linear
multi-step methods, the {non-linear} all-at-once system is
\begin{equation}\label{eq1.4}
 (B_1\otimes M){\bm u}+(B_2\times I_x)F({\bm u})={\bm b},
\end{equation}
 where $F({\bm u})=(f^\top(U_1), \dots, f^\top(U_{N_t}))^\top$.  The
 Jacobian matrix of \eqref{eq1.4} is
 \begin{equation}\label{jacobian}
  B_1\otimes M+(B_2\otimes I_x)\nabla F({\bm u}),
 \end{equation}
 where $\nabla F({\bm u})={\rm blkdiag}(\nabla f(U_1),\dots, \nabla
 f(U_{N_t}))$. To apply ParaDiag, {we approximate the} Jacobian
 matrix {\eqref{jacobian} by}
 $$
{\bm P}_\alpha({\bm u}):= C_1^{(\alpha)}\otimes M+C_2^{(\alpha)}\otimes \overline{\nabla f}({\bm u}),
 $$
 where $\overline{\nabla f}({\bm u})$ is constructed from the $N_t$
 values $\{U_n\}$ {by some} {\em averaging} \cite{GH17}, e.g.,
 $\overline{\nabla f}({\bm u})=\frac{1}{N_t}\sum_{n=1}^{N_t}\nabla
 f(U_n)$ or $\overline{\nabla f}({\bm u})=\nabla
 f(\frac{1}{N_t}\sum_{n=1}^{N_t}U_n)$.  Then, we can solve
 \eqref{eq1.4} by the following simplified Newton iteration:
 \begin{equation}\label{eq1.5}
{\bm P}_\alpha({\bm u}^{k-1})\Delta {\bm u}^{k-1}=-\left((B_1\otimes M){\bm u}^{k-1}+(B_2\times I_x)F({\bm u}^{k-1})-{\bm b}\right),~{\bm u}^{k}={\bm u}^{k-1}+\Delta {\bm u}^{k-1},
\end{equation}
where for each iteration the increment $\Delta{\bm u}^{k-1}$ can be
obtained using a ParaDiag algorithm performing the three steps in
\eqref{eq1.3}.  If we use different step sizes as in \cite{GH17},
{then} $B_1$ and $B_2$ are already diagonalizable, {and} we can
replace ${\bm P}_\alpha$ by $B_1\otimes M+(B_2\otimes
I_x)\overline{\nabla f}({\bm u}^{k-1})$ in \eqref{eq1.5}.

In practice, the ParaDiag algorithms can be combined with a {\em
  windowing} technique: after a certain number of time steps computed
in parallel in the current time window, the computation can be
restarted for the next time window in a sequential way. {This
  permits the use of a certain} adaptivity in time and space.

\section{ParaDiag for Linear  Advection-Diffusion Problems}\label{sec1}

{To illustrate the ParaDiag-I and ParaDiag-II algorithms, we now}
use the {concrete example of the} advection-diffusion equation
{with periodic boundary conditions\footnote{We use periodic
    boundary condition make the advection dominated situation harder
    for PinT algorithms, see \cite{GV08Hyp}.}}
\begin{equation}\label{eq2.1}
\begin{cases}
u_t-\nu u_{xx}+u_x=0, &(x,t)\in(-1,1)\times(0, T),\\
u(-1,t)=u(1,t), &t\in(0, T),\\
u(x,0)=e^{-30x^2}, &x\in(-1,1),
\end{cases}
\end{equation}
where $\nu>0$. {Using the} {\em method of lines} and {a}
centered finite difference scheme for the spatial derivatives, we get
the system of ODEs
\begin{subequations}
\begin{equation}\label{eq2.2a}
\dot{U}(t)+AU(t)=0,~U(0)=U_0,
\end{equation}
where the matrix $A\in\mathbb{R}^{N_x\times N_x}$ is
\begin{equation}\label{eq2.2b}
A=\frac{\nu}{\Delta x^2}\begin{bmatrix}2 &-1 &  & &-1\\
-1 &2 &-1 & & & \\
 &\ddots &\ddots &\ddots &\\
 & &-1 &2 &-1\\
-1 &  &   &-1 &2
\end{bmatrix}+
\frac{1}{2\Delta x}\begin{bmatrix}0 &1 &  & &-1\\
-1 &0 &1 & & & \\
 &\ddots &\ddots &\ddots &\\
 & &-1 &0 &1\\
1 &  &   &-1 &0
\end{bmatrix}.
\end{equation}
Here $N_x=\frac{2}{\Delta x}${, and the periodic boundary conditions
cause a zero eigenvalue in the matrix $A$.}
\end{subequations}

\subsection{ParaDiag-I}\label{sec1.1}

{ParaDiag-I consists of direct PinT solvers, and so far there two
  members {in this class of time parallel algorithms: ParaDiag-I with variable step sizes, and ParaDiag-I using hybrid time discretizations}}.

\subsubsection{ParaDiag-I: Using Variable Step Sizes}\label{sec2.1.1} 
{To use ParaDiag as a direct solver, one {can} use different time
  steps to make the time stepping matrix diagonalizable, and one
  possibility is to use} geometrically increasing time step
sizes\footnote{{Another possibility suggested by Nick Higham after
    a presentation of the first author is to use random step sizes,
    but first numerical experiments did not show an advantage over the
    geometrically increasing time steps.}}  $\{\Delta t_n\}$ to
discretize \eqref{eq2.2a} {as proposed in \cite{MR08}},
\begin{equation}\label{eq2.3}
\Delta t_n=\Delta t_1\tau^{n-1}, n\geq1,
\end{equation}
where $\tau>1$ is {a} free parameter and $\Delta t_1$ is the first step
size.  We use {as example here} the linear $\theta$-method,
\begin{equation}\label{eq2.4}
\frac{U_{n+1}-U_{n}}{\Delta t_{n+1}}+A[\theta U_{n+1}+(1-\theta)U_{n}]=0,~n=0,1,\dots, N_t-1,
\end{equation}
{and} will only consider $\theta=1$ and $\theta=\frac{1}{2}$, which
corresponds to the Backward-Euler method and the Trapezoidal rule. For
$\theta=\frac{1}{2}$, the method is also called {the}
Crank-Nicolson scheme.  The $N_t$ difference equations \eqref{eq2.4}
can be {combined} into the {\em all-at-once} system
\begin{subequations}
\begin{equation}\label{eq2.5a}
\left(B_1\otimes I_x+B_2\otimes A\right){\bm u}={\bm b},
\end{equation}
where ${\bm u}=(U_1^\top, \dots, U_{N_t}^\top)^\top$, $I_x\in\mathbb{R}^{N_x\times N_x}$ is an identity matrix and    $B_1, B_2\in\mathbb{R}^{N_t\times N_t}$  are matrices representing  the time-discretization{, namely}
\begin{equation}\label{eq2.5b}
B_1=\begin{bmatrix}\frac{1}{\Delta t_1} &  & &\\
-\frac{1}{\Delta t_2} &\frac{1}{\Delta t_2} & &\\
 &\ddots &\ddots &\\
& &-\frac{1}{\Delta t_{N_t}} &\frac{1}{\Delta t_{N_t}}
\end{bmatrix},~B_2=\begin{bmatrix}\theta   & & &\\
1-\theta &\theta & &\\
 &\ddots &\ddots &\\
& &1-\theta &\theta
\end{bmatrix}.
\end{equation}
\end{subequations}
The right hand-side ${\bm b}$ is given by
$ {\bm b}=(b_1^\top, 0,\dots, 0)^\top$ with $b_1=
\left(\frac{I_x}{\Delta t_1}-(1-\theta) A\right)U_0$.

Let $B{:=}B_2^{-1}B_1$ and $\tilde{\bm b}{:=}(B_2^{-1}\otimes
I_x){\bm b}$.  Then, we can rewrite \eqref{eq2.5a} as
\begin{equation}\label{eq2.6}
\left(B\otimes I_x+I_t\otimes A\right){\bm u}=\tilde{\bm b},
\end{equation}
where $I_t\in\mathbb{R}^{N_t\times N_t}$ is an identity matrix.  The
diagonalization of $B$ for $\theta=1$ and $\theta=\frac{1}{2}$ can be
found in \cite{GH16} and \cite{GH19} respectively{, but for the}
reader's convenience, we show the details {here:}
\begin{theorem}[see \cite{GH16,GH19}]
For the geometrically increasing step sizes $\{\Delta t_n\}$ given by
\eqref{eq2.3} with $\tau>1$, the matrix $B$ can be diagonalized as
$B=VDV^{-1}$, where $D={\rm diag}(\frac{1}{\theta\Delta_1}, \dots,
\frac{1}{\theta\Delta_{N_t}})$.  The eigenvector matrix $V$ and its
inverse are Toeplitz matrices of the form
 $$
 V=\begin{bmatrix}
 1 & & & &\\
 p_1 &1 & & &\\
  p_2 &p_1 &1 & &\\
 \vdots &\ddots &\ddots &\ddots &\\
 p_{N_t-1}    &\dots &p_2 &p_1 &1
 \end{bmatrix},~
  V^{-1}=\begin{bmatrix}
 1 & & & &\\
 q_1 &1 & & &\\
  q_2 &q_1 &1 & &\\
 \vdots &\ddots &\ddots &\ddots &\\
q_{N_t-1}    &\dots &q_2 &q_1 &1
 \end{bmatrix},
 $$
 where
 $$
 \begin{cases}
 p_n=\frac{1}{\prod_{j=1}^n(1-\tau^j)}, ~q_n=(-1)^n\tau^{\frac{n(n-1)}{2}}p_n, &\theta=1,\\
 p_n= {\prod_{j=1}^n\frac{1+\tau^j}{1-\tau^j}}, ~q_n=q^{-n}{\prod_{j=1}^n\frac{1+\tau^{-j+2}}{1-\tau^{-j}}},, &\theta=\frac{1}{2}.
 \end{cases}
 $$
\end{theorem}
Now {using the typical ParaDiag} factorization
$$
B\otimes I_x+I_t\otimes A=(V\otimes I_x)(D\otimes I_x+I_t\otimes A)(V^{-1}\otimes I_x),
$$
we can solve \eqref{eq2.6} {by performing the three steps}
\begin{equation}\label{eq2.7}
\begin{split}
&\text{Step-(a)} ~~S_1=(V^{-1}\otimes I_x)\tilde{\bm b},\\
&\text{Step-(b)} ~~S_{2, n}=\left(\frac{1}{\theta\Delta t_n}+A\right)^{-1}S_{1,n}, ~n=1,2,\dots,N_t,\\
&\text{Step-(c)} ~~{\bm u}=(V\otimes I_x)S_2,\\
\end{split}
\end{equation}
where $S_1=(S_{1,1}^\top,\dots, S_{1,N_t}^\top)^\top$ and
$S_2=(S_{2,1}^\top,\dots, S_{2,N_t}^\top)^\top$.  Since $V$ and
$V^{-1}$ are given {in closed form}, we only have to do {matrix
  vector} multiplications for Step-(a) and Step-(c){, or one could
  use a fast Toeplitz solver based on Fourier techniques}. For
Step-(b), the $N_t$ linear systems can be solved simultaneously in
parallel. {There is however an important issue with this direct
  time parallel solver ParaDiag-I: if the time steps are very
  different, the truncation error of the time stepping scheme becomes
  worse, and if they are very close to each other, ParaDiag-I suffers
  from roundoff error in the diagonalization used in} Step-(a) and
Step-(c). {The best one can do is to balance the two errors, as a
  detailed analysis in \cite{GH16,GH19} shows, and this limits the
  applicability of ParaDiag-I to shorter time intervals and few time
  steps:} the roundoff error is proportional to the condition number
of $V$, i.e.,
$$
\texttt{roundoff error}\propto\text{Cond}_2(V).
$$
If $V$ is an eigenvector matrix of $B$, the {\em scaled} matrix
$\widetilde{V}=V\widetilde{D}$ with any invertible diagonal matrix
$\widetilde{D}$ is an eigenvector matrix of $B$ as well. From
\cite{GH16,GH19}, the matrix $\widetilde{D}={\rm
  diag}\left((1+\sum_{j=1}^{N_t-n}|p_j|^2)^{-\frac{1}{2}}\right)$ is a
good choice.

\begin{figure}[ht]
  \centering
  \includegraphics[width=0.45\textwidth,height=0.3\textwidth]{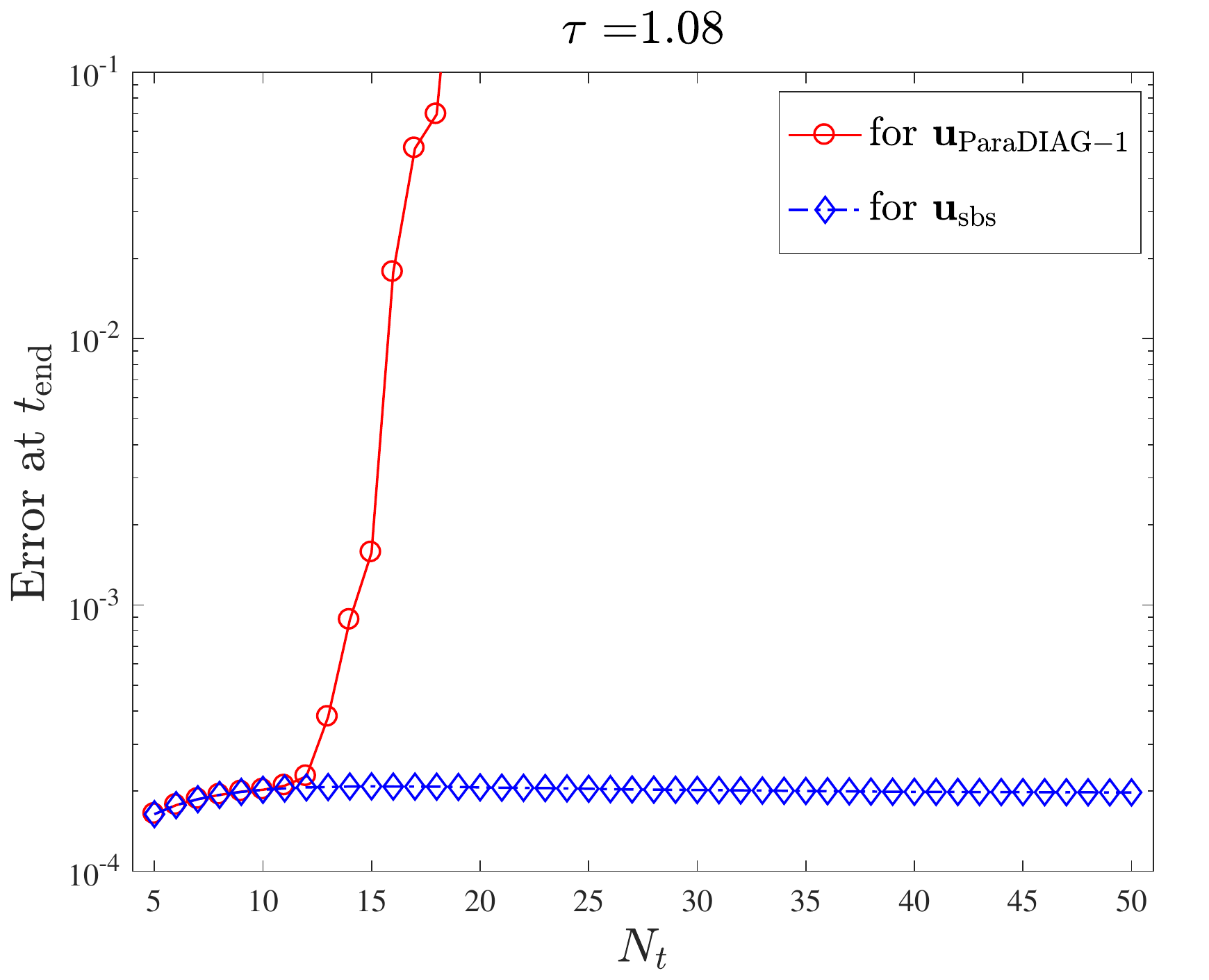}
  \includegraphics[width=0.45\textwidth,height=0.3\textwidth]{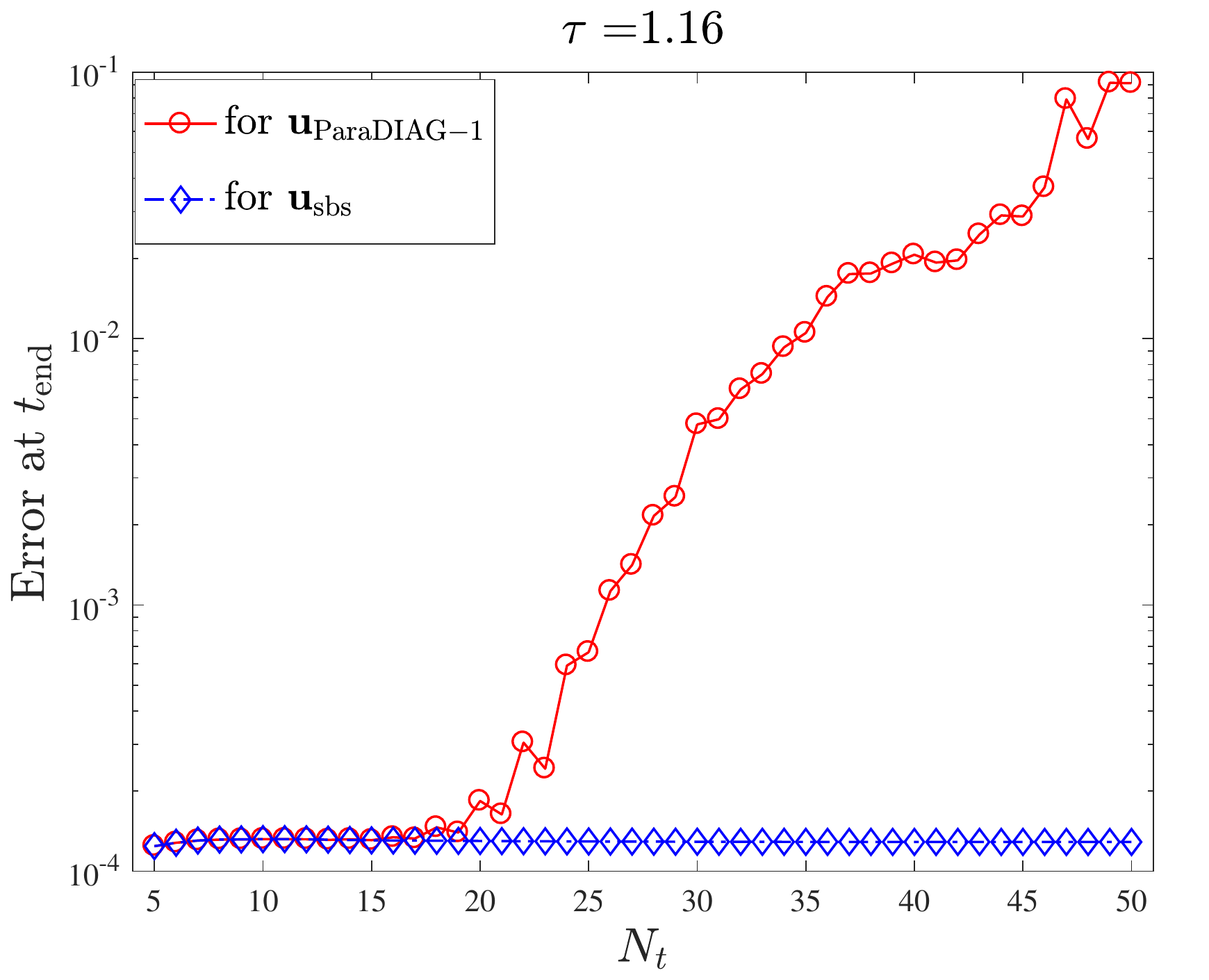}
  \caption{\em Using the geometric time steps \eqref{eq2.8} with
     $\Delta t_{N_t}=10^{-2}$  and two different values of $\tau$, the errors measured at the finial time point
    $t=t_{\rm end}$  for two numerical solutions: ${\bm u}_{\rm sbs}$ obtained
    step by step (dash-dot lines) and ${\bm u}_{\rm ParaDiag-I}$
    obtained by \eqref{eq2.7} (solid lines).  Here, $\nu=10^{-2}$ and the Trapezoidal rule is used.
  } \label{fig1.1}
\end{figure}

\begin{figure}[ht!]
  \centering
 \includegraphics[width=0.32\textwidth]{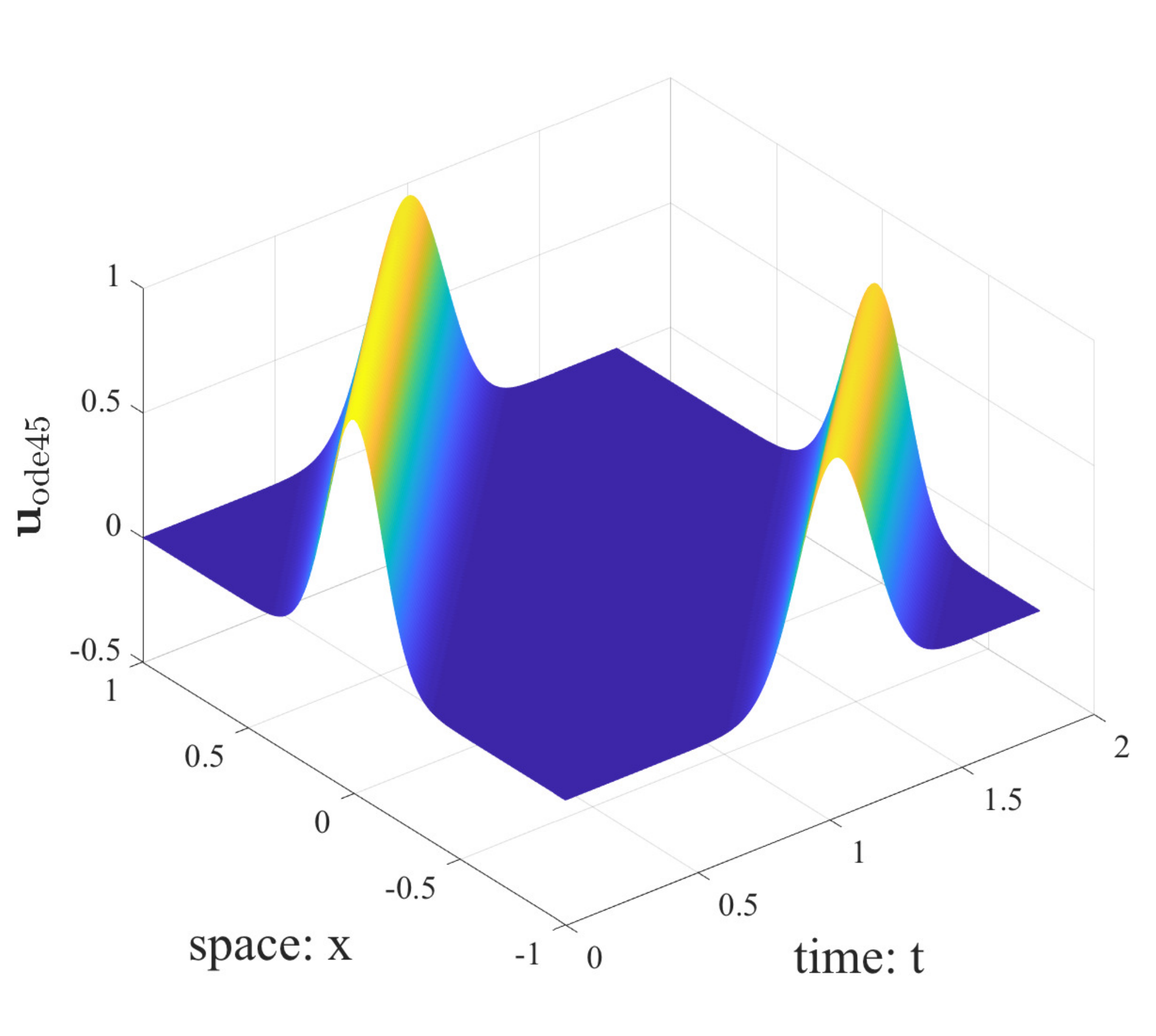} \includegraphics[width=0.32\textwidth]{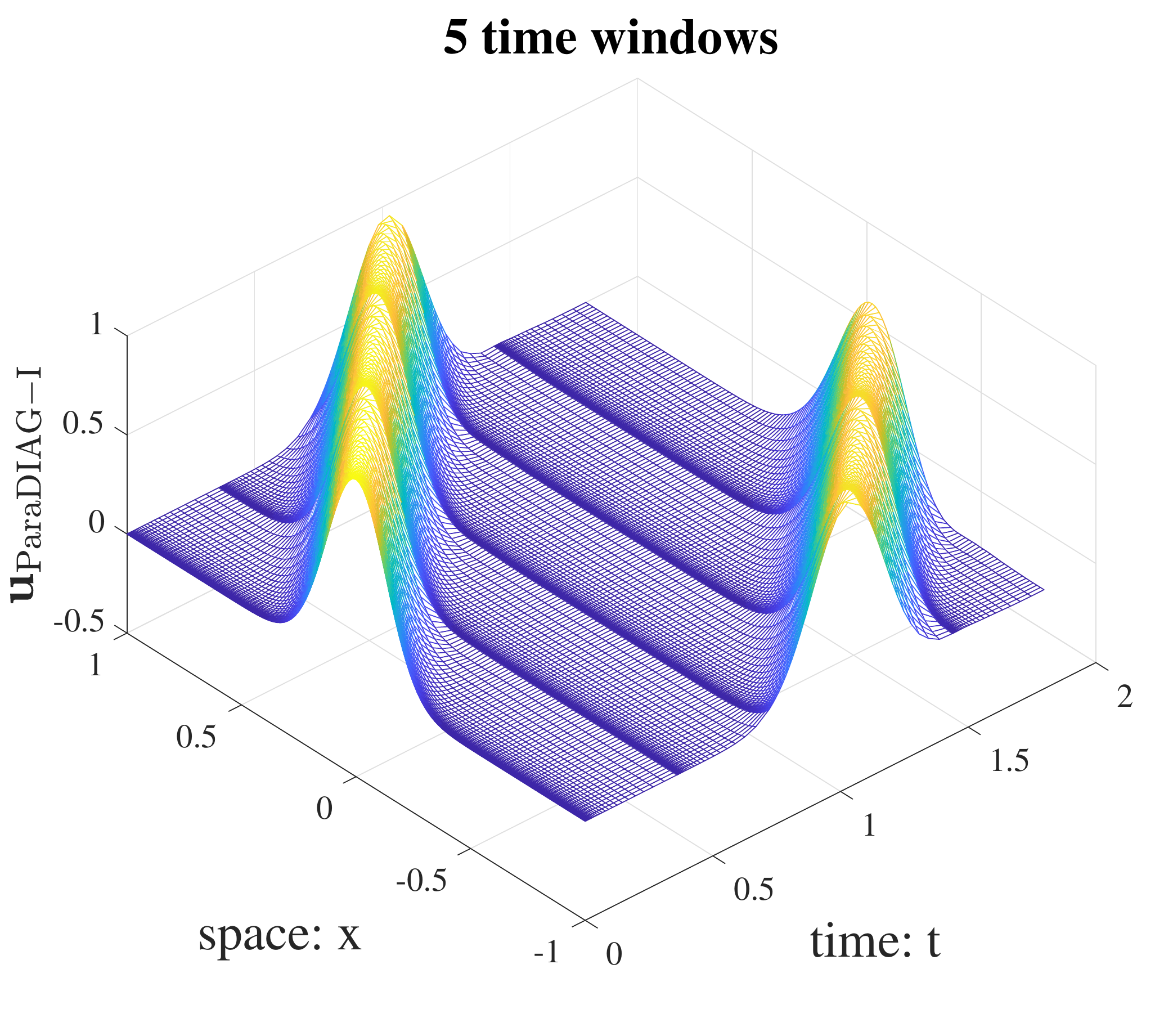} \includegraphics[width=0.32\textwidth]{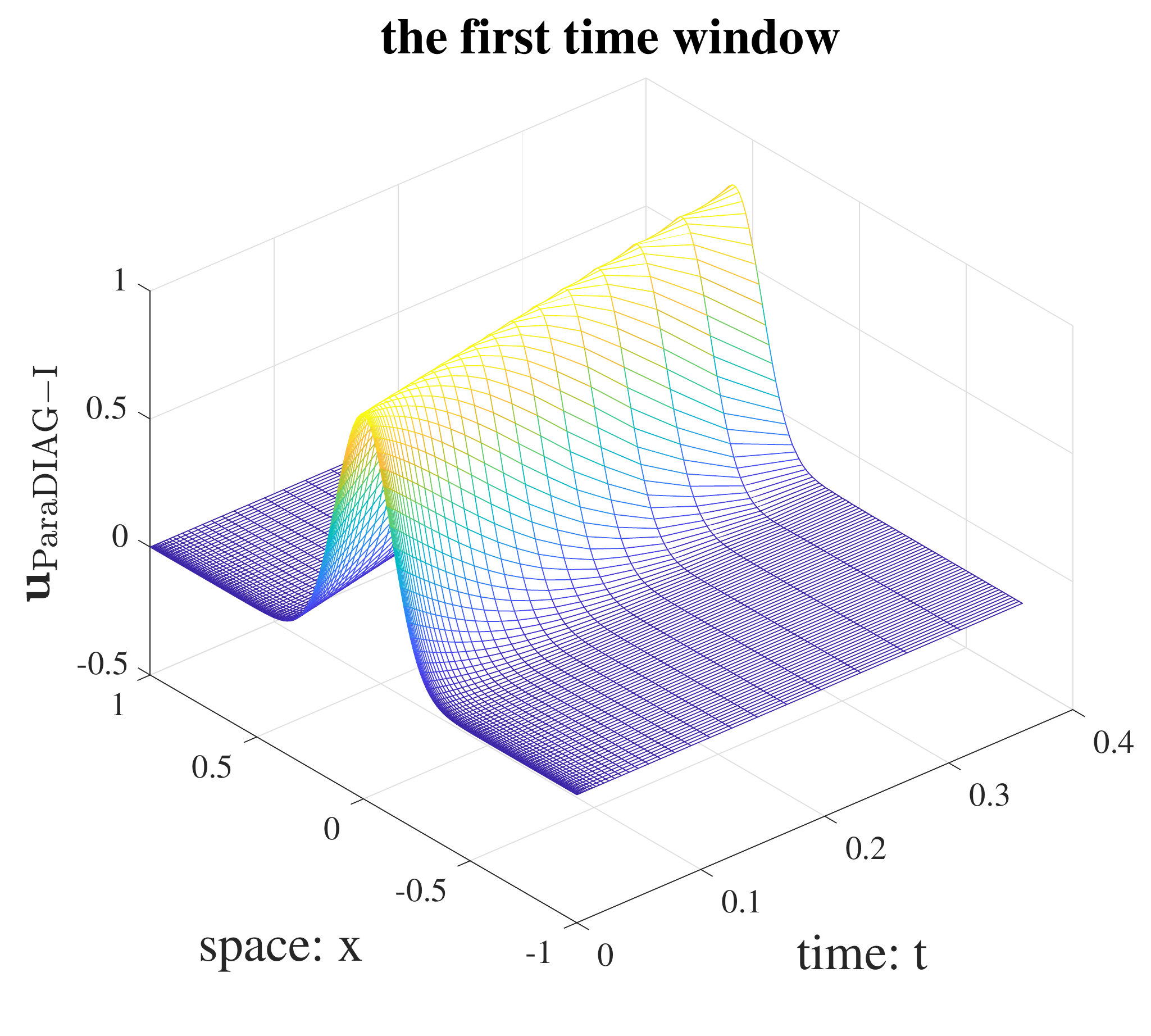}
  \caption{\em Left: reference solution ${\bm u}_{\rm ode45}$.
    Middle: numerical solution ${\bm u}_{\rm ParaDiag-I}$ in 5 time
    windows. Right: the numerical solution ${\bm u}_{\rm ParaDiag-I}$
    in the first time window, which clearly shows the graded mesh
    sizes in time.  Here, $\nu=10^{-3}$ and the Trapezoidal rule is
    used. For ${\bm u}_{\rm ParaDiag-I}$, each time window contains
    $N_t=32$ steps and the parameters $\tau$ and $\Delta t_{N_t}$ in
    \eqref{eq2.8} are $\tau=1.16$ and $\Delta t_{N_t}=0.05$.
  } \label{fig1.2}
\end{figure}

To illustrate the limitations of {this first ParaDiag-I variant},
we provide the Matlab code \texttt{ParaDiag\_V1\_for\_ADE.m}, to test
it for the advection-diffusion equation.  For given $N_t$ and $\tau$
and the final step-size $\Delta t_{N_t}$ (e.g., $\Delta
t_{N_t}=10^{-2}$)|this $\Delta t_{N_t}$ determines the maximal
discretization error, we specify the first $N_t-1$ step-size $\{\Delta
t_n\}$ as
\begin{equation}\label{eq2.8}
\Delta t_n=\Delta t_{N_t}\times \tau^{n-N_t}, ~n=1,2,\dots, N_t-1.
\end{equation}
For the space discretization, we fix $\Delta x=\frac{1}{64}$. To study
the accuracy of {this ParaDiag-I variant,} we use a reference
solution ${\bm u}_{\rm ode45}$ {obtained from the Matlab ODE solver}
\texttt{ode45} with a very small absolute and relative tolerance,
\texttt{AbsTol}=$10^{-12}$ and \texttt{RelTol}=$10^{-12}$.  In Figure
\ref{fig1.1}, we show the measured error at the end time point $t_{\rm
  end}$ for $\mathbf{u}_{\rm sbs}$ and $\mathbf{u}_{\rm ParaDiag-1}$
as $N_t$ increases.
(For given
$\tau$ and $\Delta t_{N_t}$, such a $t_{\rm end}$ grows as $N_t$
increases.)  We clearly see that using the geometric time steps
\eqref{eq2.8} {degrades the} accuracy of the numerical solution, and
when the time steps are too similar, the roundoff error problem sets
in. This phenomenon was carefully studied in \cite{GH16,GH19}, and the
best possible geometrically stretched grid was determined, which leads
to precise limits of time window length and number of time steps
within which {this original ParaDiag-I variant} can be reliably
used.
The accuracy shown in Figure \ref{fig1.1} indicates that
the number of time steps, i.e., $N_t$, can not be large.  For long
time computation, we can divide the whole time interval into several
time windows and each time window includes a moderate number {of} time
steps, say $N_t=20\sim 30$. Then, we apply ParaDiag-I to these time
windows one by one.  An illustration of such a {\em windowing}
technique is shown in Figure \ref{fig1.2}.

\subsubsection{ParaDiag-I: Using a Hybrid Time Discretization (\red{\bf New Progress)}}\label{sec2.1.2}

{We have seen} that using variable {time} step sizes {poses}
a serious restriction on the number of time {steps} $N_t$ {one
  can use}. In practice, the algorithm only works for $N_t=20\sim 30$
{for our model problem}. {We present now a new} direct PinT
solver {in the ParaDiag-I family}, for which much larger
{numbers of time steps} $N_t$ can be used.

With a uniform step-size $\Delta t$, we use an explicit {mid}-point
scheme for the first $(N_t-1)$ time steps followed by an implicit
Euler method for the last step, that is
\begin{subequations}
\begin{equation}\label{neq2.9a}
\begin{cases}
\frac{U_{n+1}-U_{n-1}}{2\Delta t}+AU_n=0, ~n=1,2,\dots, N_t-1,\\
\frac{U_{N_t}-U_{N_t-1}}{\Delta t}+AU_{N_t}=0.
\end{cases}
\end{equation}
The all-at-once {system} of this scheme is 
 \begin{equation}\label{neq2.9b}
{\bm A}{\bm u}={\bm b},~{\bm A}=B\otimes I_x+I_t\otimes A, 
\end{equation}
where 
\begin{equation}\label{neq2.9c}
B=\frac{1}{\Delta t}\begin{bmatrix}
0 &\frac{1}{2} & & & \\
-\frac{1}{2} &0 &\frac{1}{2}  & & \\ 
  &\ddots &\ddots &\ddots &\\ 
 &  &-\frac{1}{2}  &0   &\frac{1}{2}\\
 &  &   &-1   &1\\
\end{bmatrix},~ {\bm b}=\begin{bmatrix}\frac{u_0}{2\Delta t}\\ 0\\ \vdots\\ 0\end{bmatrix},~{\bm u}=\begin{bmatrix}U_1\\ U_2\\ \vdots\\ U_{N_t}\end{bmatrix}. 
\end{equation}
\end{subequations}  

Such a hybrid time discretization is not new and a brief history is as
follows.  In 1985, Axelsson and Verwer \cite{AV85} studied this scheme
with the aim of circumventing the well-known Dahlquist-barriers
between convergence and stability which arise in using the
{mid}-point scheme {for} time-stepping.  In the general
nonlinear case they proved that the numerical solutions obtained
simultaneously are of uniform second-order accuracy (see Theorem 4 in
\cite{AV85}), even though the last step is {only first-order}.
Numerical results in \cite{AV85} indicate that the hybrid time
discretization \eqref{neq2.9a} is suitable for stiff problems in both
linear and nonlinear cases. A very similar hybrid time discretization
{was} investigated {already} by Fox in 1954 \cite{F54} and Fox
and Mitchell in 1957 \cite{FM57}, where instead of backward-Euler
the authors use {a} BDF2 method for the last step:
$$
 \frac{3U_{N_t}-4U_{N_t-1}+U_{N_t-2}}{2\Delta t}+AU_{N_t}=0.
$$
In this case, the time discretization matrix $B$ however {looses}
the tridiagonal structure and {certain} desirable properties to
efficiently handle the all-at-once system do not hold any more.  This
{subject} was carefully {studied} by Brugnano, Mazzia and
Trigiante in 1993 \cite{BMT93}, who focus on solving the all-at-once
system \eqref{neq2.9b} iteratively {with a}
preconditioner {performing} two operations: a block odd-even cyclic
reduction of ${\bm A}$ and a scaling procedure for the result{ing} matrix
by its diagonal blocks. The block cyclic reduction requires
matrix-matrix multiplications {involving} $A$, and the scaling requires
to invert $I_x+4\Delta t^2A$ and $I_x+2\Delta tA(I_x+\Delta tA)$. Both
operations are expensive if $A$ arises from semi-discretiz{ations of} a PDE in
high dimension and/or with fine mesh sizes.  Nowadays, the hybrid time
discretization \eqref{neq2.9a} is a famous example of the so-called
{\em boundary value methods} (BVMs) \cite{BT98}.

\begin{figure}
\centering
\includegraphics[width=0.45\textwidth]{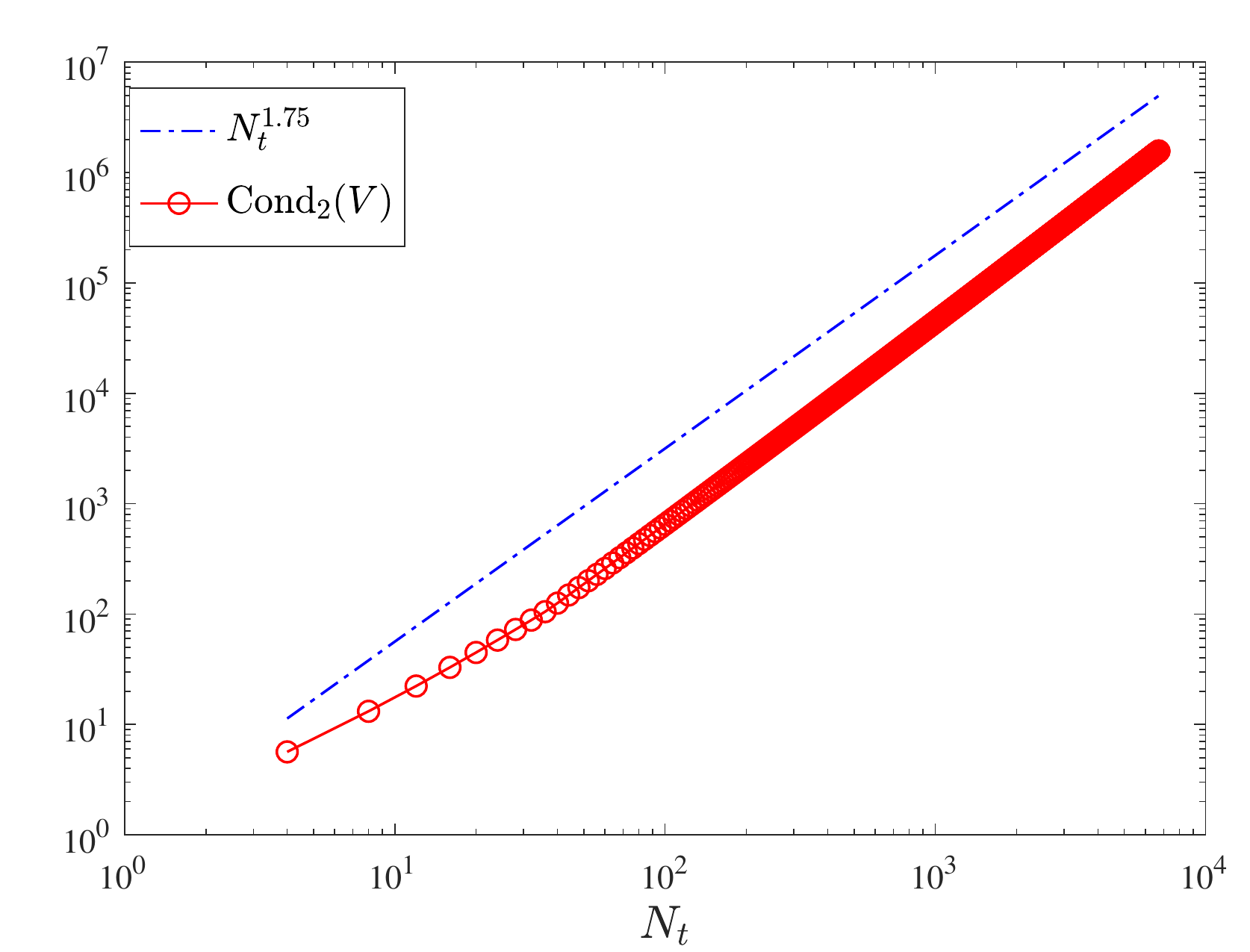} \includegraphics[width=0.45\textwidth]{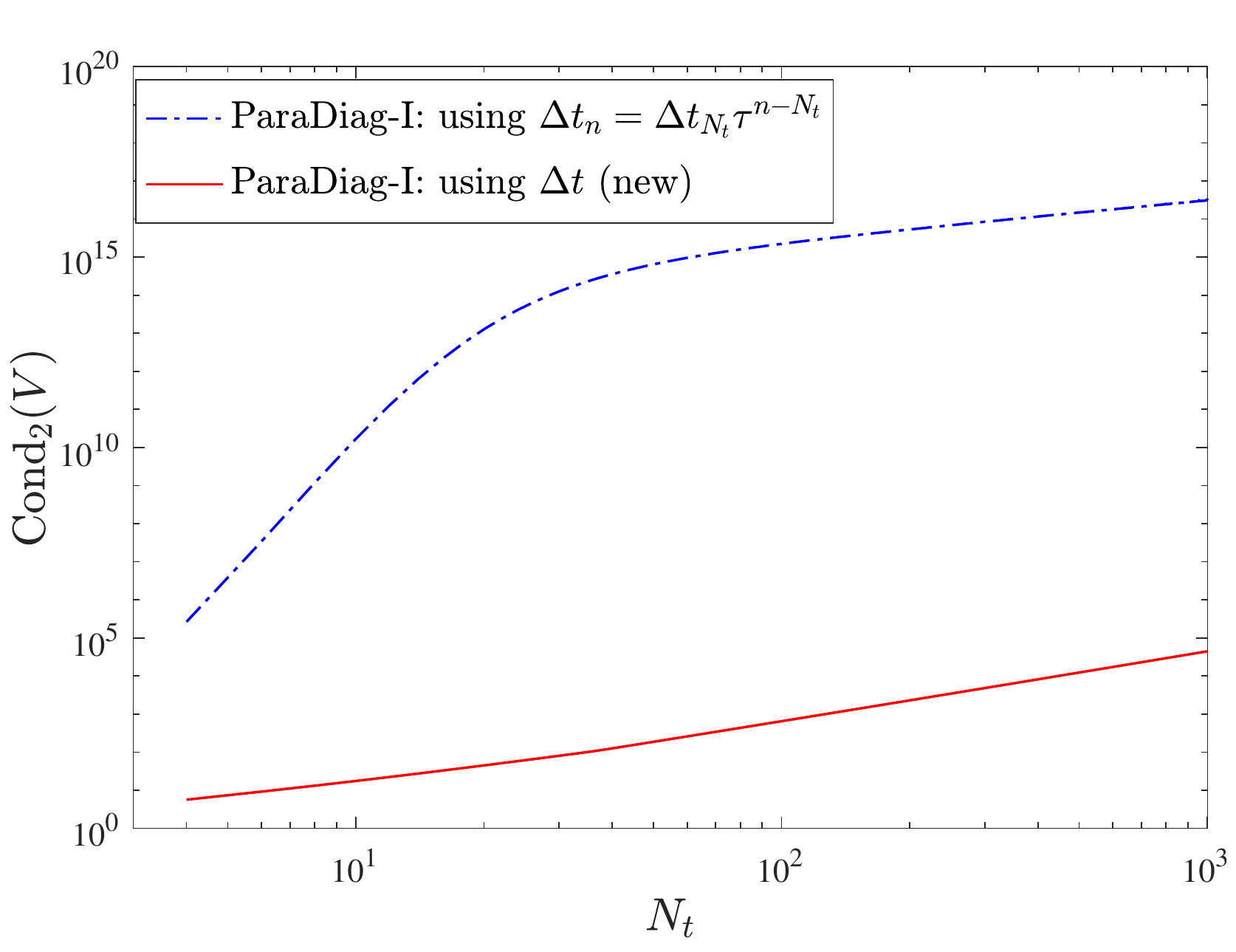} 
\caption{\em In practice, the condition number of the eigenvector matrix $V$ satisfies ${\rm Cond}_2(V)=\mathcal{O}(N_t^{1.75})$, which is better than the theoretical bound \eqref{condUorder}. } \label{nfig2.3}
\end{figure}

The reason for using \eqref{neq2.9a} as the time integrator for ParaDiag-I is the following interesting property. 
\begin{theorem}[\cite{LWW21}] \label{thmCondU}
  Let ${\rm i}=\sqrt{-1}$ be the imaginary unit and {let}
  $$
  P^{(1)}_{N_t}(x)=\cos(N_t\arccos x),~
  P^{(2)}_{N_t}(x)=\sin[(N_t+1)\arccos x]/\sin(\arccos x),
  $$
  be the Chebyshev polynomials of first and second kind
  of degree $N_t$.  The matrix $B$ in \eqref{neq2.9c} can be
  diagonalized as $B=VDV^{-1}$ with eigenvalues {$\lambda_n={\rm
      i}\Delta t x_n$, where $\{x_n\}_{n=1}^{N_t}$} are the $N_t$
  different roots of
  \begin{equation}\label{e-value}
    P^{(2)}_{N_t-1}(x)-{\rm i}  P^{(1)}_{N_t}(x)=0.
  \end{equation} 
  For $\lambda_n$, the corresponding eigenvector $\bm v_n=[v_{n,0},
    \cdots,v_{n,n-1}]^\T$ is given {by}
   \begin{equation}\label{e-vector}
 		v_{n,l}={\rm i}^lP^{(2)}_l(x_n),~ l=0,\cdots,N_t-1,
   \end{equation}
   where $v_{n,0}=1$ is assumed for normalization.  For $N_t\ge 8$, 
   {the condition number satisfies}
   \begin{equation} \label{condUorder}
	{\rm Cond}_2(V)= \mathcal{O}(N_t^2).
   \end{equation} 
\end{theorem}
In practice, the eigenvalues $\{\lambda_n\}$ can be obtained by
applying a Newton iteration to \eqref{e-value} and then the
eigenvector matrix $V$ {is given by}
\begin{equation}  \label{Vformula}
 V=[\bm v_1,\bm v_2, \cdots, \bm v_{N_t}]={\rm diag}\left({\rm i}^0,{\rm i}^1,\cdots,{\rm i}^{N_t-1}\right)
\begin{bmatrix}  
 	P^{(2)}_0(x_1)&\cdots & P^{(2)}_0(x_{N_t}) \\
   	\vdots &\cdots   &\vdots \\
 	P^{(2)}_{N_t-1}(x_1) &\cdots & P^{(2)}_{N_t-1}(x_{N_t}) \\
 	\end{bmatrix}.
\end{equation}  
So, Theorem \ref{thmCondU} actually provides a close{d} form
eigendecomposition of the time stepping matrix $B$ in \eqref{neq2.9c}.
 The estimate
Cond$_2(V)=\mathcal{O}(N_t^2)$ is {on the conservative side,} and
in practice we find ${\rm Cond}_2(V)=\mathcal{O}(N_t^{1.75})$; see
Figure \ref{nfig2.3} on the left for an illustration{, but we
  currently do not have a prove of this.}
For comparison, in Figure
\ref{nfig2.3} on the right we show Cond$_2(V)$ for the two direct PinT
algorithms, the one {studied} in \cite{GH19} {going back to
  \cite{MR08}}, and the new one introduced in \cite{LWW21}.

\begin{figure}[h]
\centering
\includegraphics[width=0.45\textwidth]{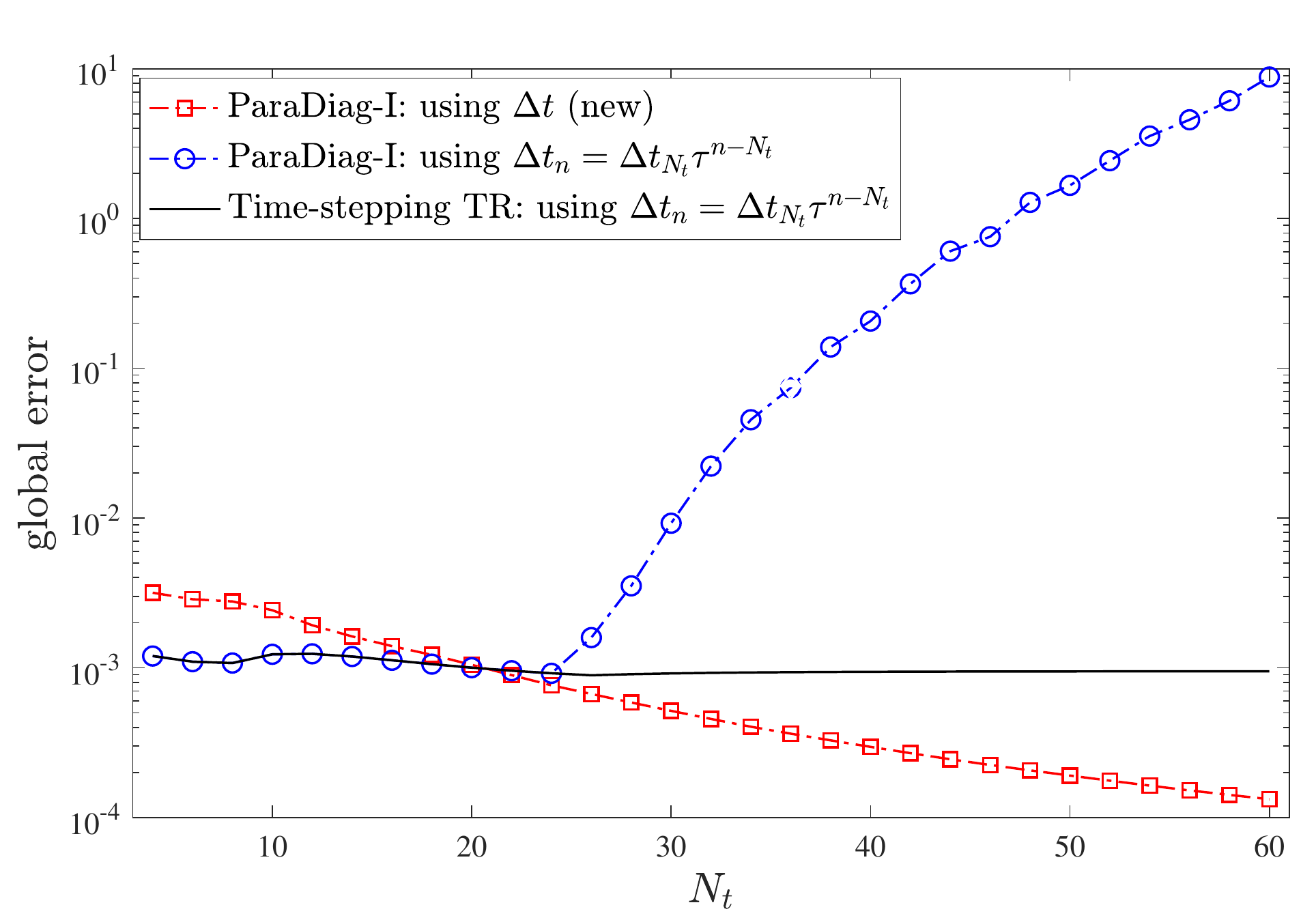}
\caption{\em The global error for the new ParaDiag-I algorithm
  introduced in \cite{LWW21}, the ParaDiag-I algorithm in \cite{GH19}
  and the time-stepping TR using the variable step-sizes.
} \label{nfig2.4}
\end{figure}

The we{a}k dependence of the condition number on $N_t$ implies that the
roundoff error arising from the diagonalization procedure only
moderately increases as $N_t$ grows.  To validate this, we consider
the advection-diffusion equation with $\nu=1e-3$. The ODE system
\eqref{eq2.2a} is obtained by semi-discretizing \eqref{eq2.1} with
$\Delta x=\frac{1}{128}$.  For the ParaDiag-I algorithm in Section
\ref{sec2.1.1}, we use the Trapezoidal rule (TR) as the
time-integrator, where the step-sizes are fixed {to} $\Delta t_n=\Delta
t_{N_t}\tau^{n-N_t}$ with $\Delta t_{N_t}=10^{-2}$ and $\tau=1.15$.
For each $N_t$ we {run} the ParaDiag-I algorithm in \cite{GH19} by
using the variable step-sizes. Then, we calculate the length of the
time interval, i.e., $T(\tau, N_t)=\sum_{n=1}^{N_t}\Delta
t_j$\footnote[1]{For the algorithm in \cite{GH19}, since $\Delta
  t_n=\Delta t_{N_t}\tau^{n-N_t}$ the length of time interval grows as
  $N_t$ increases.  } and {run} the new algorithm by using a
uniform step-size $\Delta t={T(\tau, N_t)}/{N_t}$. {We} define the global
error as
$$
\mbox{global error}=\max_{n=1,2,\dots, N_t}\|{\bm u}_{n,h}-{\bm u}^{\rm ref}_{n,h}\|_\infty, 
$$
where $\{{\bm u}^{\rm ref}_{n,h}\}$ denotes the reference solution
obtained by using the \texttt{expm} function in MATLAB, i.e., ${\bm
  u}^{\rm ref}_{n,h}=\texttt{expm}(-t_nA)$.  The sequence $\{{\bm
  u}_{n,h}\}$ is obtained {in three different} ways: by the new algorithm
introduced here, by the algorithm in \cite{GH19} and by the
time-stepping TR using the variable step-sizes.  The comparison for
the global error of these three numerical solutions is shown in Figure
\ref{nfig2.4}. It is clear that the global error of the new algorithm
continuously decrease{s} when $N_t\geq6$. {This confirms that} the
condition number shown in Figure \ref{nfig2.3} {leads to a}
roundoff error much
smaller than the time discretization error, and thus the global error
is dominated by the time discretization error.  The decreas{e} of the
global error is because the {time} step-size
$$
\Delta t= \Delta t_{N_t}{\sum}_{n=1}^{N_t}\tau^{n-N_t}/{N_t}=\Delta t_{N_t}\frac{1-\tau^{-N_t}}{n(1-\tau^{-1})}\approx\frac{0.0766}{N_t} ~(\text{if}~N_t\geq40)
$$
 decreases as $N_t$ grows and thus the time discretization error decreases as well.   

For    second-order  problems 
\begin{subequations}
\begin{equation}\label{neq2.14a}
U''+AU=F, ~U(0)=U_0, U'(0)=\bar{U}_0, 
\end{equation}
 we can first  represent {them} as a first-order system (with $V=U'$)
\begin{equation}\label{neq2.14b}
 W'+\begin{bmatrix} & -I_x\\ A &\end{bmatrix}W=\begin{bmatrix}0\\ g\end{bmatrix}, ~W:=\begin{bmatrix}U\\ V\end{bmatrix}, 
\end{equation}
\end{subequations}
and then apply the new ParaDiag-I algorithm to \eqref{neq2.14b}.
However, this doubles the system size and storage requirement for the
space variables at each time point.  To avoid this, we eliminate the
auxiliary variables $\{V_n\}$ at the discrete level and only solve the
all-at-once system for $\{U_n\}$. Let $B_{\rm 2nd}$ be the
corresponding time discretization matrix {which has the interesting
  property}
\begin{equation}\label{neq2.15}
 B_{\rm 2nd}=B^2= \frac{1}{4\Delta t^2}\begin{bmatrix}
	-1 &0 &1 & & & & \\
	0 &-2 &0  &1 & & & \\
	1 &0 &-2 &0    &1 & & \\
	&\ddots &\ddots &\ddots &\ddots  &\ddots \\ 
	& &1 &0 &-2 &0 &1\\
	& & &1 &0&-3 &2\\
	& & & &2&-4 &2\\
\end{bmatrix}{,} 
\end{equation}
where $B$ is the matrix in \eqref{neq2.9c} (the time stepping matrix
for the first-order problem \eqref{neq2.9a}).  Thus, the same
eigendecomposition of $B$ can be reused and the condition number of
the eigenvector matrix is not effected.
 
{We show next} a parallel experiment for {a} 2D wave equation defined on
$\Omega=(0, 1)^2$,
\begin{equation}\label{wavePDE}
	\begin{cases}
		u_{tt}(x,y,t)- \Delta u(x,y,t) =f(x,y,t), &\tn{in} \Omega\times(0,T),\\
		u(x,y,t)=0,&\tn{on} \partial\Omega\times(0, T), \\
		u(x,y,0)=u_0(x,y),&\tn{in} \Omega,\\
		u_t(x,y,0)=\bar u_0(x,y),&\tn{in} \Omega,
	\end{cases}
\end{equation} 
with the {initial conditions and right hand side}
\begin{equation*}
\begin{split}
&u_0(x,y)=0, ~ \bar u_0(x,y)=2\pi  x(x-1)y(y-1),\\
&f(x,y,t)=-4\pi^2x(x-1)y(y-1)\sin(2\pi t)-2\sin(2\pi t)(x(x-1)+y(y-1)).
\end{split}
\end{equation*}
The exact solution of this problem is $u(x,y,t)=x(x-1)y(y-1)\sin(2\pi
t).$ Using the centered finite difference scheme, we obtain a
second-order ODE system \eqref{neq2.14a} with $A\approx-\Delta$ being
the spatial discrete matrix. We show in Table \ref{Ex2Tab1} the
approximation errors {and} strong and weak scaling results for the new
ParaDiag-I algorithm (the CPU time is measured in seconds).
\begin{table}[ht] 
	\centering
	\caption{Scaling Results {for} a Wave {equation} ($T=2$ with $\Delta x=\Delta y=1/512$)}
	
	\begin{tabular}{|c||c|c||c|c|c||c|c|c|c||c|c|c|c|c|cc}
		\hline
		Core\#&\multicolumn{5}{c|}{Strong scaling} & \multicolumn{4}{c|}{Weak scaling} \\ \hline
		$s$&	$N_t$&  Error& CPU &Sp.&SE & $N_t$&  Error& CPU&WE\\
		\hline  
		1       & 512   & 7.88e-05     &1328.6   &1.0     &   100.0\%               &   2   & 9.19e-03        & 5.4 &  100.0\%   \\ 
		2        & 512   & 7.88e-05    &676.3    &2.0     &   98.2\%               &   4   & 2.21e-02        & 5.4  &  100.0\%  \\ 
		4        & 512   & 7.88e-05    &332.6    &4.0     &   99.9\%               &   8   & 3.16e-01        & 5.5  &  100.0\% \\ 
		8        & 512   & 7.88e-05    &172.6    &7.7     &   96.2\%               &  16   & 1.33e-01        & 5.7  &  100.0\% \\ 
		16      & 512   & 7.88e-05     &91.2     &14.6    &   91.0\%                &  32   & 2.30e-02        & 6.0  & 94.8\%  \\ 
		32      & 512   & 7.88e-05     &51.7     &25.7    &   80.3\%                &  64   & 5.21e-03        & 7.1  & 82.1\%  \\ 
		64      & 512   & 7.88e-05     &31.2     &42.6    &   66.5\%                & 128   & 1.27e-03        & 9.5  & 67.9\%  \\ 
		128     & 512   & 7.88e-05     &23.2     &57.3    &   44.7\%                & 256   & 3.16e-04        & 14.8 & 46.6\%  \\ 
		256     & 512   & 7.88e-05     &20.3     &65.4    &   25.6\%                & 512   & 7.88e-05        & 27.4 & 28.2\%  \\
		\hline
	\end{tabular}
	\label{Ex2Tab1}
\end{table}
{These results were obtained on the
SIUE Campus Cluster} with 10 CPU nodes connected via {a} 25-Gigabit per
second (Gbps) Ethernet network, where each node is equipped with two
AMD EPYC 7F52 16-Core Processors at 3.5GHz base clock and 256GB
RAM. {The notations} in Table \ref{Ex2Tab1} are \[ \mbox{Speedup
  (Sp.)}=\frac{\texttt{CT}(N_t,1)}{\texttt{CT}(N_t,s)}{,}
\]
where \texttt{CT}$(N_t,s)$ {is} the measured CPU (wall-clock) time by using $s$ cores for $N_t$ time {steps}. 
 The strong and weak scaling efficiency with $s$ cores {ís} computed as
\[
\mbox{Strong Efficiency (SE)}= \frac{\texttt{CT}(N_t,1)}{s\times \texttt{CT}(N_t,s)}
,\qquad
\mbox{Weak Efficiency (WE)}=\frac{\texttt{CT}(2,1)}{\texttt{CT}(2 \times s,s)}.
\]
  In Table
\ref{Ex2Tab1}, both {the} strong and weak scaling efficienc{ies}
are very promising up to 32 cores{, but} when the core number
$s\geq64$, we see a drop of the parallel efficiency. This is mainly
due to the slow interconnect between the nodes (each node contains 32
cores).  We remark that the measured parallel speedup and efficiency
are affected by many factors, such as computer cluster setting and how
the parallel codes {are implemented}. {Our} parallel results
{here} may {still} underestimate the best possible speedup and
efficiency with optimized parallel codes.
{Our Matlab} and parallel codes
(\texttt{ParaDIAG\_V1\_Hybrid\_for\_ADE.m},
\texttt{ParaDIAG\_V1\_Hybrid\_for\_Wave.c}) used here can be found {at}
\url{https://github.com/wushulin/ParaDIAG}.

\subsection{ParaDiag-II}

Instead of using ParaDiag as a direct solver, we can use it
iteratively and solve a nearby problem in each iteration chosen such
that the time stepping matrix of the nearby problem ({even} with uniform
time step size) can still be diagonalized. This idea leads to ParaDiag
algorithms in the ParaDiag-II group.
  {In} this
group, we can use ParaDiag within a stationary iteration, or as a
{preconditioner for a} Krylov subspace method.  There are so far
two very different ways to use ParaDiag within a stationary iteration,
proposed in \cite{GW19} and \cite{W18}. The use of ParaDiag {as a
  preconditioner for} a Krylov subspace method can be found in
\cite{MPW18,WL20b}.

\subsubsection{{ParaDiag-II -- Waveform Relaxation (WR) Variant}}

The ParaDiag algorithm introduced in \cite{GW19} is based on the
{Waveform Relaxation} iteration
 \begin{equation}\label{eq2.9}
\dot{U}^k(t)+AU^k(t)=0,~U^k(0)=U_0+\alpha(U^k(T)-U^{k-1}(T)),~ t\in(0, T),
\end{equation}
where $k\geq1$ is the iteration index and $\alpha\in(0, 1]$ is a free
  parameter.  Upon convergence, the tail term
  $\alpha(U^k(T)-U^{k-1}(T))$ is canceled and thus the converged
  solution is the solution of \eqref{eq2.2a}.  Applying the linear
  $\theta$-method with a uniform step size $\Delta t$ to \eqref{eq2.9}
  gives
\begin{equation}\label{eq2.10}
\begin{cases}
\frac{U^k_{n}-U^k_{n-1}}{\Delta t}+A\left(\theta U^{k}_n +(1-\theta)U^k_{n-1}\right)=0,~n=1,2,\dots, N_t,\\
U^k_0=\alpha U^k_{N_t}-\alpha U^{k-1}_{N_t}+U_0,
\end{cases}
\end{equation}
where $N_t=T/\Delta t$. We rewrite \eqref{eq2.10} as an {\em
  all-at-once} system,
 \begin{subequations}
\begin{equation}\label{eq2.11a}
\begin{split}
\left(C_1^{(\alpha)}\otimes I_x+C_2^{(\alpha)}\otimes A\right){\bm u}^k={\bm b}^{k-1},
\end{split}
\end{equation}
{with} ${\bm u}^k=(U^k_1, \dots, U^k_{N_t})^\top$, {and}
${C}_{1}^{(\alpha)}, {C}_2^{(\alpha)}\in\mathbb{R}^{N_t\times N_t}$
and ${\bm b}^{k-1}\in\mathbb{R}^{N_t N_x}$ are given by
\begin{equation}\label{eq2.11b}
\begin{split}
&C_1^{(\alpha)}=
\frac{1}{\Delta t}\begin{bmatrix}
1 & & &-\alpha\\
-1 &1 & &\\
&\ddots &\ddots &\\
& &-1 &1
\end{bmatrix},~
C_2^{(\alpha)}=\begin{bmatrix}
\theta & & &(1-\theta)\alpha\\
1-\theta &\theta & &\\
&\ddots &\ddots &\\
& &1-\theta &\theta
\end{bmatrix},\\
&{\bm b}^{k-1}=\left((U_0-\alpha U^{k-1}_{N_t})\left(\frac{1}{\Delta t}I_x-(1-\theta)A\right),0,\dots,0\right)^\top.
\end{split}
\end{equation}
\end{subequations}
The matrices $C^{(\alpha)}_{1,2}$ are so-called $\alpha$-circulant
matri{ces} and can be diagonalized as stated in Lemma
\ref{lem1}{, and we can again use the typical ParaDiag
  factorization} $C_1^{(\alpha)}\otimes I_x+C_2^{(\alpha)}\otimes A=
(V\otimes I_x)\left(D_1\otimes I_x+D_2\otimes A\right)(V^{-1}\otimes
I_x)$. Hence, similar to \eqref{eq2.7} we can solve \eqref{eq2.11a}
{performing the three steps}
\begin{equation}\label{eq2.12}
\begin{split}
&\text{Step-(a)} ~~S_1=(\mathbb{F}\otimes I_x)(\Gamma_\alpha\otimes I_x){\bm b}^{k-1},\\
&\text{Step-(b)} ~~S_{2, n}= (\lambda_{1,n}I_x+\lambda_{2,n}A)^{-1}S_{1,n}, ~n=1,2,\dots,N_t,\\
&\text{Step-(c)} ~~{\bm u}^k=(\Gamma_\alpha^{-1}\otimes I_x)(\mathbb{F}^*\otimes I_x)S_2,\\
\end{split}
\end{equation}
{where $D_j=\text{diag}(\lambda_{j,1},\dots,\lambda_{j,N_t})$ and
  $j=1,2$.}  In (\ref{eq2.7}), Step-(a) and Step-(c) can be computed
efficiently via FFT and Step-(b) is {again} highly parallel.  The
eigenvector matrix $V$ satisfies
\begin{equation}\label{eq2.13}
{\rm Cond}_2(V)={\rm Cond}_2(\Gamma_\alpha^{-1}\mathbb{F}^*)\leq {\rm Cond}_2(\Gamma_\alpha^{-1}){\rm Cond}_2(\mathbb{F}^*)={\rm Cond}_2(\Gamma_\alpha^{-1})\leq \frac{1}{\alpha},
\end{equation}
{and thus the conditioning is depending on the choice of $\alpha$.
  The convergence properties of this ParaDiag-II algorithm are
  summarized in the following theorem.}
\begin{theorem}[{see} \cite{GW19}]\label{the1.2}
  For {the linear system of} ODEs $\dot{U}(t)+AU(t)=f$, suppose
  $\Re(\lambda(A))\geq r\geq0$ with $\lambda(A)$ being an arbitrary
  eigenvalue of $A$.  Let ${\bm u}^k$ be the $k$-th iterate of the
  ParaDiag-II algorithm \eqref{eq2.10} with $\alpha\in(0, 1)$ and
  ${\bm u}$ be the reference solution obtained by directly applying
  the same time-integrator to the {system of} ODEs. Then {the
    linear convergence estimate} $\|{\bm u}^k-{\bm u}\|_\infty\leq
  \rho^k\|{\bm u}^0-{\bm u}\|_\infty$ {holds}, where
\begin{equation*}
\rho\leq
\begin{cases}
\frac{\alpha e^{-Tr}}{1-\alpha e^{-Tr}}, &\text{\rm Backward-Euler},\\
\frac{\alpha}{1-\alpha}, &\text{\rm Trapezoidal~rule}.
\end{cases}
\end{equation*}
\end{theorem}
{This shows that the ParaDiag-II algorithm \eqref{eq2.10} converges
  with a rate} independent of the spectrum of the matrix $A$ and the
step size of {the} time-discretization.  The convergence factor
$\rho$ {becomes smaller when} $\alpha$ decreases, but the condition
number of $V$ (cf. \eqref{eq2.13}) implies that $\alpha$ can not be
arbitrarily small (e.g., {not of the size} $\alpha=10^{-13}$),
because in this case the roundoff error will pollute the accuracy.
The best parameter $\alpha_{\rm opt}$ is {again} the value
balancing the roundoff error and the discretization error{, like
  for the direct solver ParaDiag-I}, see \cite{GW19} for more
discussions. In practice, $\alpha=10^{-2}$ and $\alpha=10^{-3}$ are
good choices.

\begin{figure}
\centering
\mbox{\includegraphics[width=0.33\textwidth]{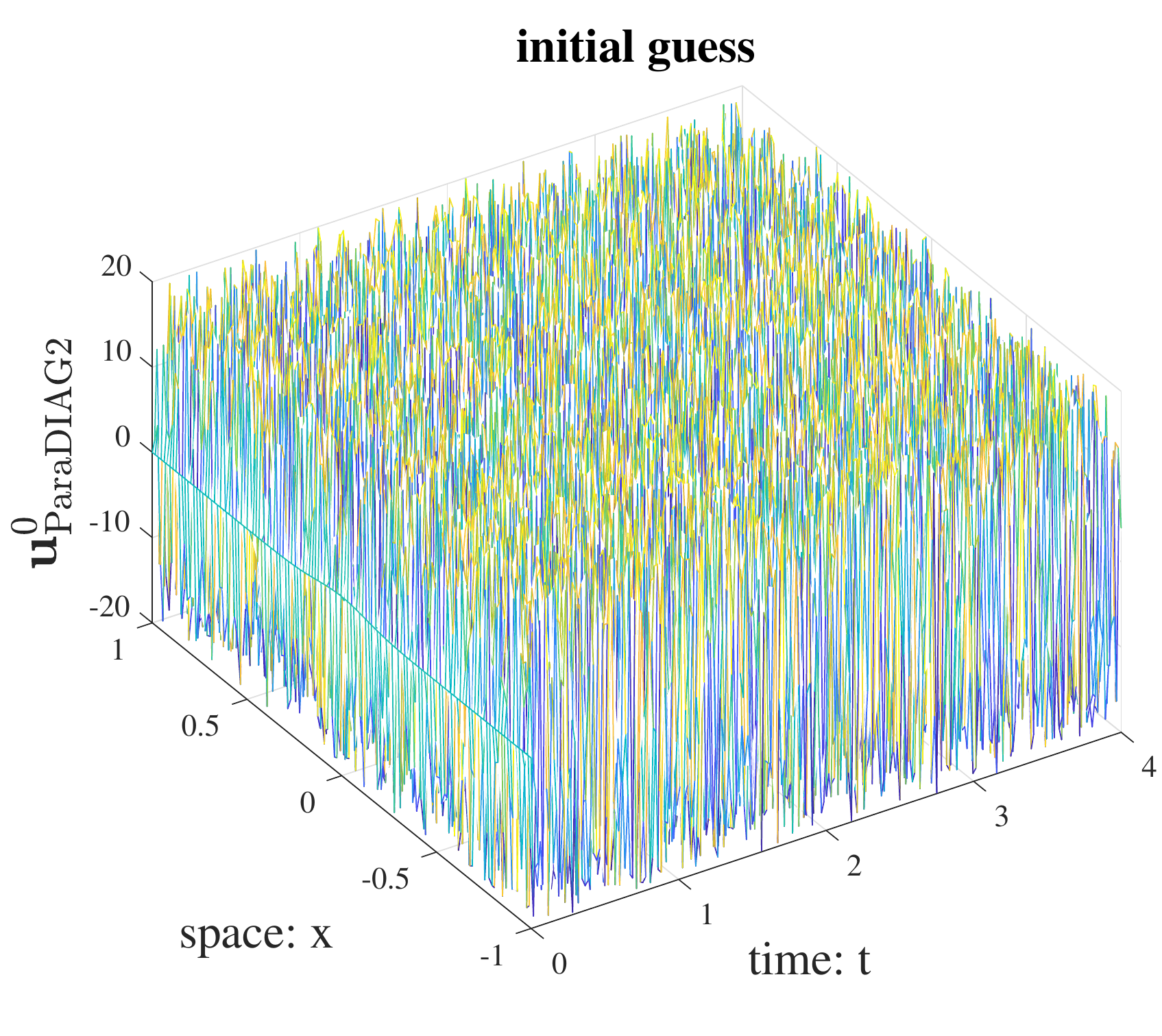}
\includegraphics[width=0.33\textwidth]{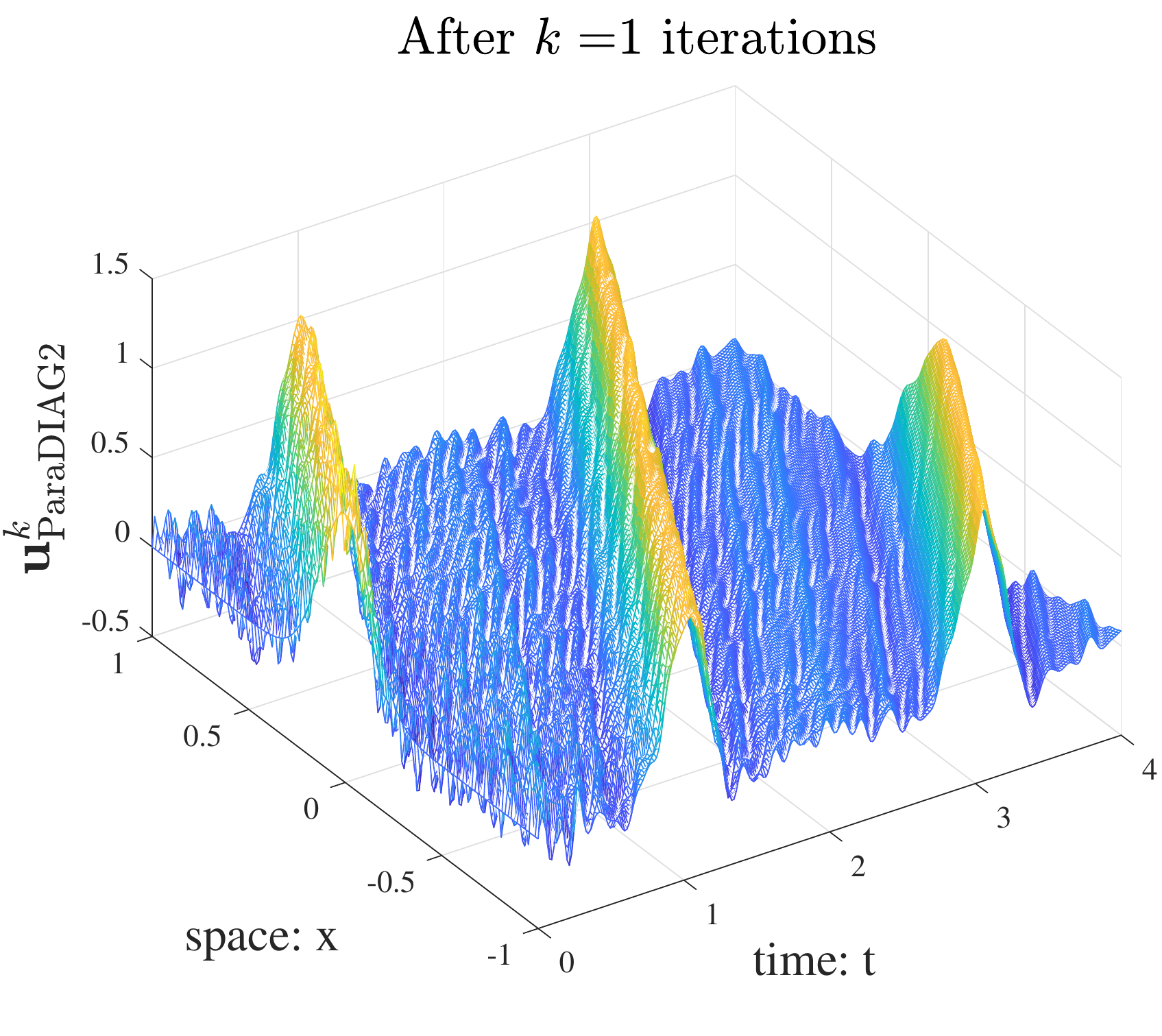}
\includegraphics[width=0.33\textwidth]{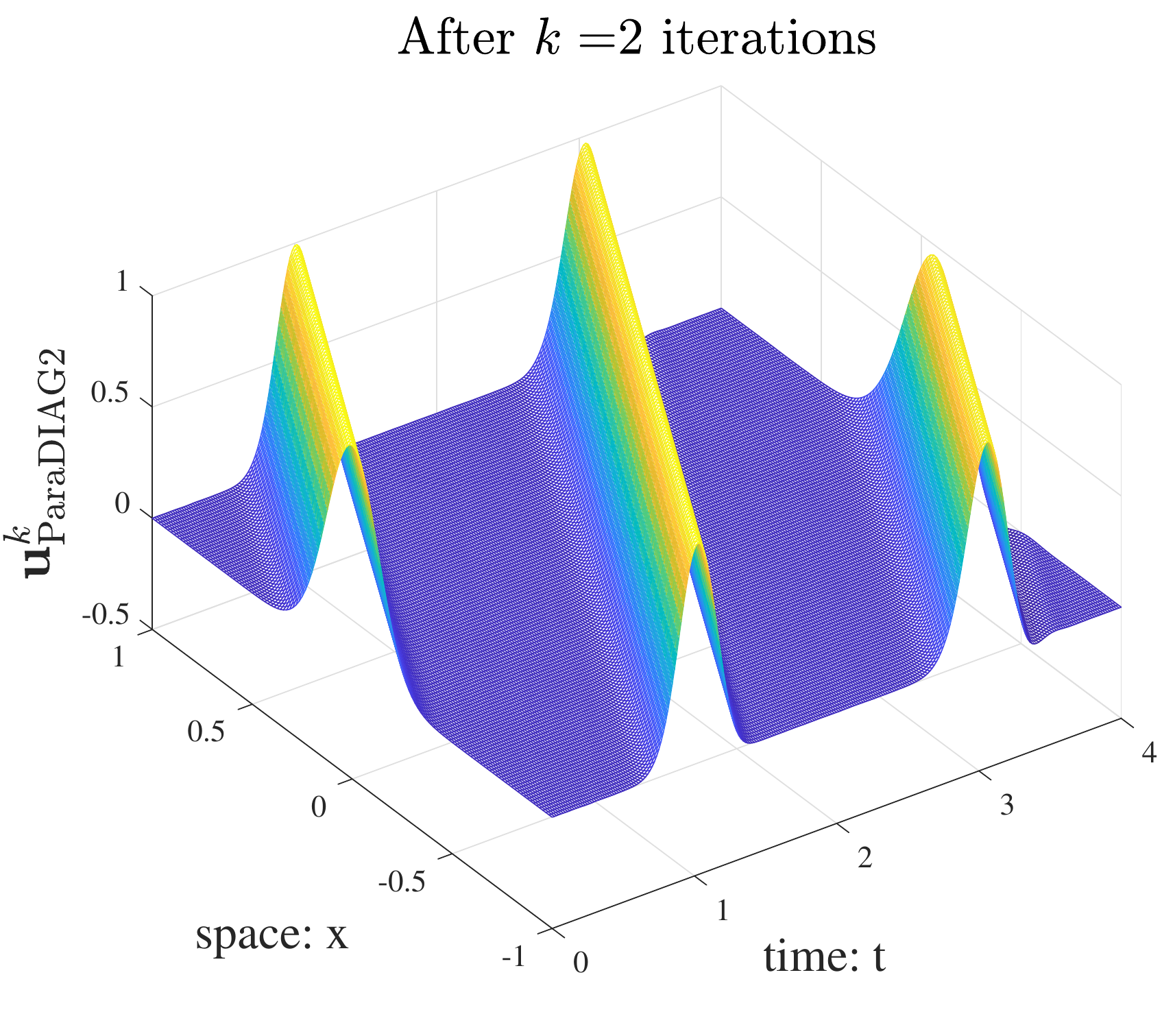}}
\caption{\em {Initial guess and the first two iterates generated by
    the ParaDiag{-II} algorithm \eqref{eq2.10} for $\nu=10^{-4}$
    and $\Delta x=\Delta t=\frac{1}{64}$}.} \label{fig1.3a}
\end{figure}

\begin{figure}[ht]
\centering
\includegraphics[width=0.49\textwidth]{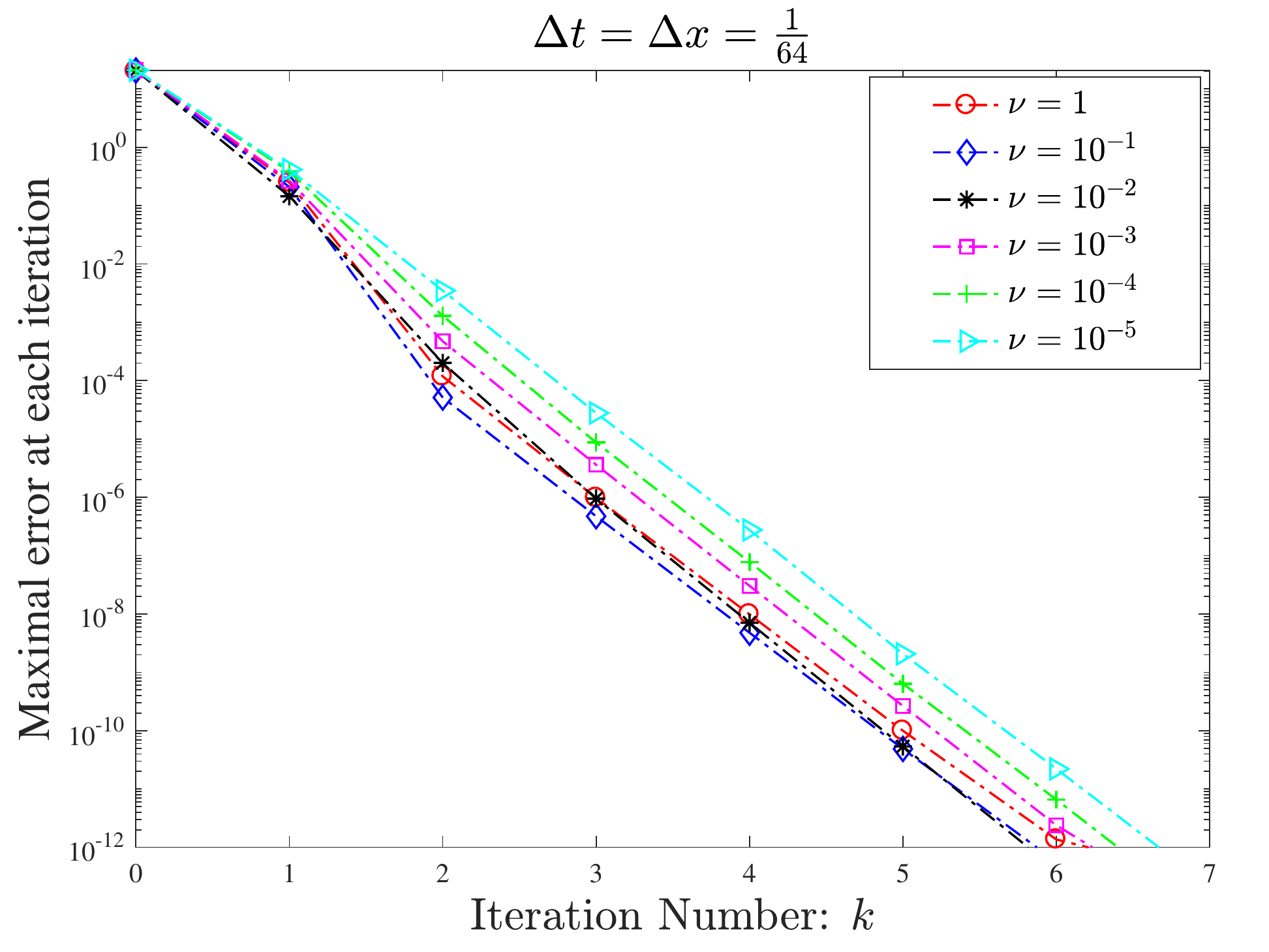}
\includegraphics[width=0.49\textwidth]{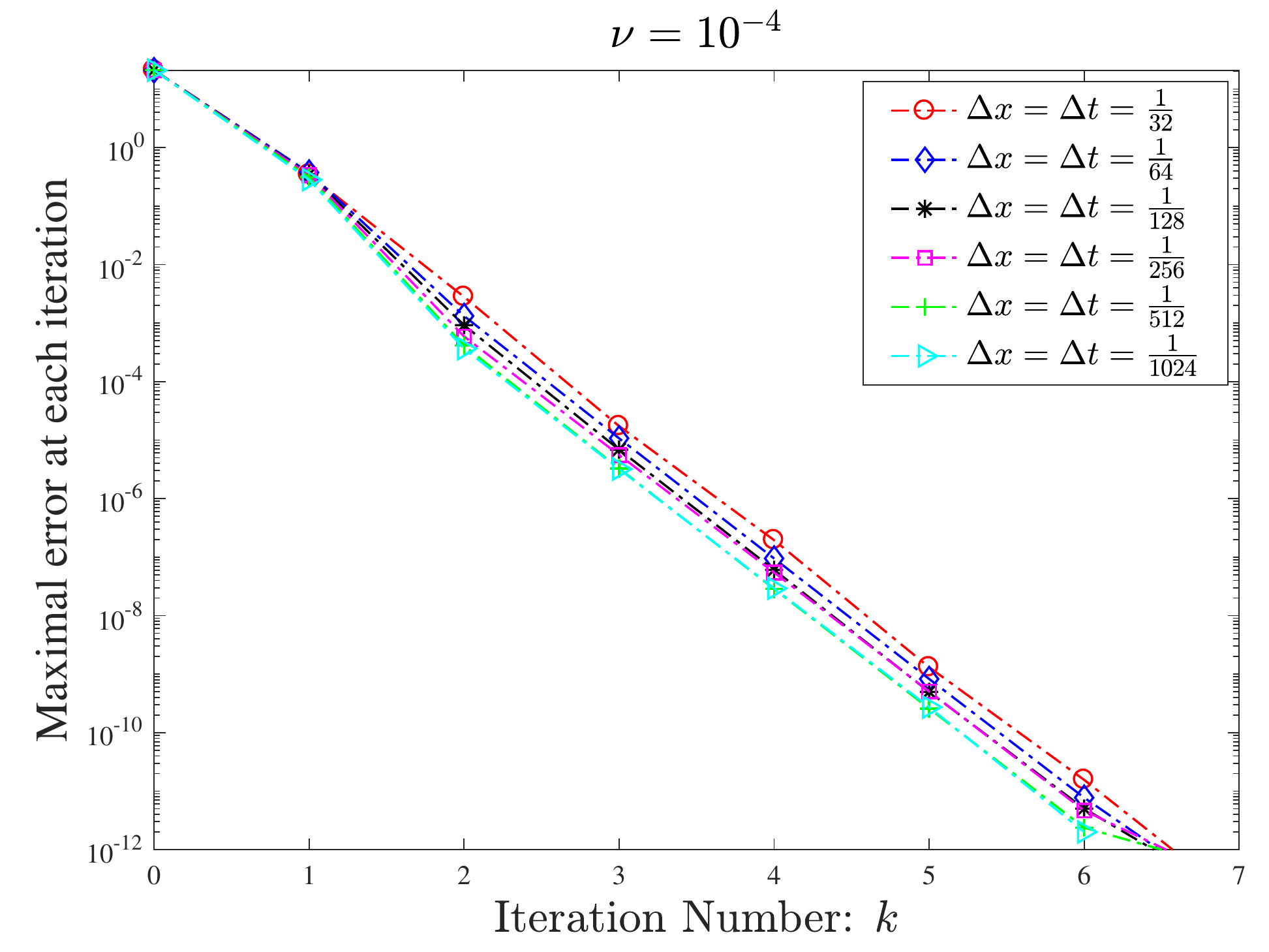}	
\caption{\em The  convergence of the ParaDiag{-II} algorithm
  \eqref{eq2.10} is robust with respect to $\nu$, $\Delta x$ and
  $\Delta t$. Here, $\alpha=10^{-2}$ and the Trapezoidal rule is used
  as the time-integrator. } \label{fig1.3b}
\end{figure}

We provide a Matlab code, namely \texttt{ParaDiag\_V2\_WR\_for\_ADE}, to
test the ParaDiag{-II} algorithm \eqref{eq2.10}. In the code, we
use the \texttt{fft} command to obtain $D_{1,2}$ by just {using}
the first columns of $C^{(\alpha)}_{1,2}$, instead of the {entire}
matrices. To implement Step-(a) in \eqref{eq2.12}, we use the \texttt{fft}
command as follows:
$$\text{\small
    b=\texttt{reshape}(b,Nx,Nt); ~sol\_stepA=\texttt{fft}(Gam.*(b.')).';
    }
$$
where b is the vector ${\bm b}^{k-1}$. Similarly, to implement
Step-(c) we use {the inverse FFT command} \texttt{ifft},
$$
\text{\small
Uk=(invGam.*\texttt{ifft}(sol\_stepB.')).';
    }
$$
Here,
Gam$={\big{(}}1,\alpha^{\frac{1}{N_t}},\dots,\alpha^{\frac{N_t-1}{N_t}}{\big{)}}$
and
invGam$={\big{(}}1,\alpha^{-\frac{1}{N_t}},\dots,\alpha^{\frac{1-N_t}{N_t}}{\big{)}}$.
With an initial guess chosen randomly as
\texttt{random}({\color{magenta}'unif'},-20,20,$N_x$, $N_t$), the
first 2 iterates {of this ParaDiag-II algorithm} are shown in
Figure \ref{fig1.3a}.
The maxim{um} error at each iteration is shown in Figure
\ref{fig1.3b}.

{We next present some parallel speedup results for the ParaDiag-II
  algorithm \eqref{eq2.10} based on Waveform Relaxation for}
a time-dependent {advection}-diffusion problem with periodic
boundary conditions in 2D,
\begin{equation}\label{eq2.14}
\begin{cases}
       \partial_t u({\bf x}, t) - \nu\Delta u({\bf x}, t) + \nabla u({\bf x}, t) = 0,
    &\textrm{in} ~ (0,T)\times \Omega, \\
        u({\bf x}, 0) = u_0({\bf x}), &{\rm in}~  \Omega,
\end{cases}
\end{equation}
where $\Omega=(0,1)\times(0,1)$ and $u_0({\bf
  x})=e^{-20[(x-\frac{1}{2})^2+(y -\frac{1}{2})^2]}$. {The results
  were} obtained on the China Tianhe-1 supercomputer \cite{YL10},
which is a multi-array, configurable and cooperative parallel system
with a theoretical peak {performance} of 1.372 petaflops, composed
of high performance general-purpose microprocessors and a high-speed
Infiniband network.  We {used} the parallel Fortran library
{MUMPS (MUltifrontal Massively Parallel sparse direct Solver
\cite{AD01,AB19})} version 4.10.0 to solve the linear systems in
Step-(b) of \eqref{eq2.12}.  For Step-(a) and Step-(c), the
\texttt{fft} and \texttt{ifft} commands are dissected into complex
arithmetic operations.  For the particular case when the source term
is zero as shown in \eqref{eq2.11b}, Step-(a) can be implemented in an
economical way: only the first column of $(\mathbb{F}\otimes
I_x)(\Gamma_\alpha\otimes I_x)$ is needed to compute $S_1$.

We {provide our} parallel codes in Fortran, which are zipped {in
  the file} 'Parallel Codes.zip'{, including} a README file, in
which we briefly introduce how to use these codes.  In addition, a
number of {comment} statements are {also contained in the}
Fortran functions and subroutines.  All the Fortran codes {were}
compiled with mpich-3.1.3 using the icc compiler version 11.1.059 and
-O2 optimization level.

\begin{figure}[ht]
\centering
\includegraphics[width=6.13in,height=4in]{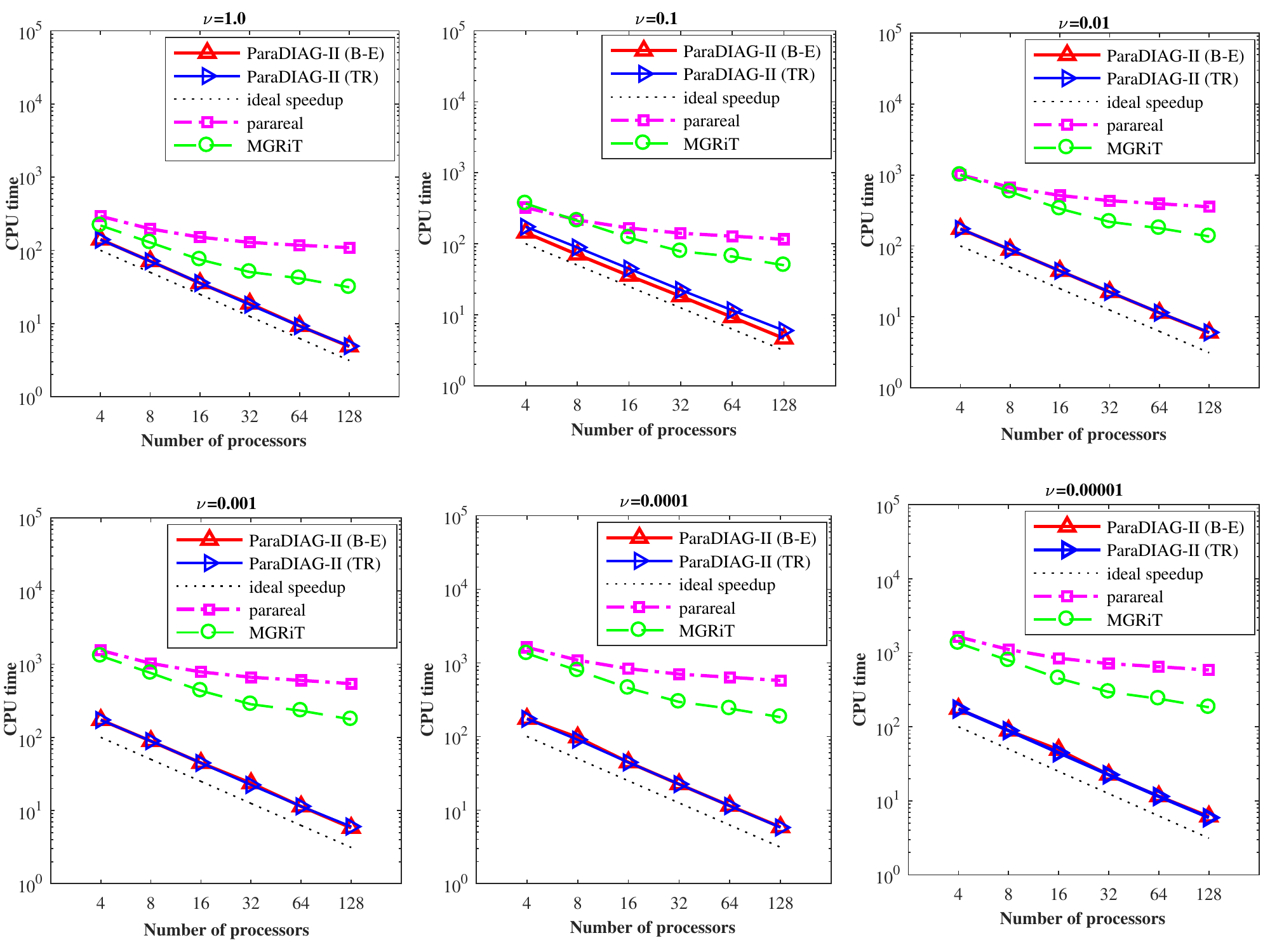}
\caption{\em Comparison of the overall time-to-solution in a strong
  scaling study, where $\Delta x=\Delta y=\Delta t=1/128$ and
  $N_t=512$.  The coarsening factor is cf=8 for both parareal and
  MGRiT.}\label{fig-04}
\end{figure}

We denote by ParaDiag-II (B-E) the algorithm \eqref{eq2.10} using
Backward-Euler, and by ParaDiag-II (TR) the one using the Trapezoidal
rule, and set $\alpha=0.02$.  For comparison, we also apply the
parareal algorithm and MGRiT to \eqref{eq2.14}.  The parareal
algorithm is implemented {using} the two-level XBraid solver with
F-relaxation, and MGRiT is the multilevel XBraid solver with
  FCF-relaxation (i.e., an initial F-relaxation followed by a
C-relaxation and then a second F-relaxation).  Furthermore, we skip
the unnecessary work during the first XBraid down cycle for both the
parareal and MGRiT
{algorithms, and fix the coarsening factor to 8}.  As shown in Table
\ref{ctp-01},
\begin{table}
  \centering
  \footnotesize
  \tabcolsep=2pt
  \caption{\em Iteration {numbers} of ParaDiag-II (B-E),
    ParaDiag-II (TR), Parareal and MGRiT in a strong scaling study,
    where $\Delta x=\Delta y=\Delta t = 1/128$ and $N_t=512$. The
    symbol $np$ indicates the number of processors and the coarsening
    factor is 8 in both parareal (two time levels) and MGRiT (three
    time levels).}\label{ctp-01}
\begin{tabular}{cccccccccccccccccccccccccccccccc}\hline
\multirow{2}{*}{\emph{np}}&\multicolumn{4}{c}{$\nu=10^0$}
&&\multicolumn{4}{c}{$\nu=10^{-1}$}&&\multicolumn{4}{c}{$\nu=10^{-2}$}
&&\multicolumn{4}{c}{$\nu=10^{-3}$}&&\multicolumn{4}{c}{$\nu=10^{-4}$}
&&\multicolumn{4}{c}{$\nu=10^{-5}$} \\ \cline{2-5} \cline{7-10} \cline{12-15} \cline{17-20} \cline{22-25} \cline{27-30}
~& B-E & TR & PR & MG && B-E & TR & PR & MG && B-E & TR & PR & MG && B-E & TR & PR & MG && B-E & TR & PR & MG && B-E & TR & PR & MG \\ \hline
   4 & 4 & 4 & 9 & 4 && 4 & 5 & 10 & 7 && 5 & 5 & 33 & 20 && 5 & 5 & 51 & 26 && 5 & 5 & 54 & 27 && 5 & 5 & 55 & 27 \\ 
   8 & 4 & 4 & 9 & 4 && 4 & 5 & 10 & 7 && 5 & 5 & 33 & 20 && 5 & 5 & 51 & 26 && 5 & 5 & 54 & 27 && 5 & 5 & 55 & 27 \\ 
  16 & 4 & 4 & 9 & 4 && 4 & 5 & 10 & 7 && 5 & 5 & 33 & 20 && 5 & 5 & 51 & 26 && 5 & 5 & 54 & 27 && 5 & 5 & 55 & 27 \\ 
  32 & 4 & 4 & 9 & 4 && 4 & 5 & 10 & 7 && 5 & 5 & 33 & 20 && 5 & 5 & 51 & 26 && 5 & 5 & 54 & 27 && 5 & 5 & 55 & 27 \\ 
  64 & 4 & 4 & 9 & 4 && 4 & 5 & 10 & 7 && 5 & 5 & 33 & 20 && 5 & 5 & 51 & 26 && 5 & 5 & 54 & 27 && 5 & 5 & 55 & 27 \\ 
 128 & 4 & 4 & 9 & 4 && 4 & 5 & 10 & 7 && 5 & 5 & 33 & 20 && 5 & 5 & 51 & 26 && 5 & 5 & 54 & 27 && 5 & 5 & 55 & 27 \\ \hline
\end{tabular}\\ \vskip 0.1cm
B-E: ParaDiag-II (B-E), TR: ParaDiag-II (TR), PR: parareal, MG: MGRiT
\end{table}
ParaDiag-II (B-E), ParaDiag-II (TR), parareal and MGRiT converge
robustly with respect to the number of processors.  ParaDiag-II (B-E)
and ParaDiag-II (TR) lead to parameter-robust convergence, while for
parareal and MGRiT the required iteration counts increase dramatically
as $\nu$ changes from 1 to $10^{-4}$.  The tolerance \texttt{tol} for
all experiments here is set to $10^{-6}(<\min\frac{\{\Delta t^2,
  \Delta x^2\}}{10}$). {In} Figure \ref{fig-04}
we compare the measured CPU time{s} for these PinT algorithms.
Clearly, ParaDiag-II (B-E) and ParaDiag-II (TR) are two {optimally}
scaling PinT algorithms, while for MGRiT and parareal the scaling is a
little bit worse (this is because of {the} sequential
coarse-grid-correction as we will {see} in the next subsection).
The {corresponding} data is given in Table \ref{ctp-02}.
\begin{table}[t]
  \centering
  \footnotesize
  \tabcolsep=4pt
  \caption{\em {Results for the} strong scaling speedup of the
    four PinT algorithms measured by $T_{\rm cpu}^{(4)} / T_{\rm
      cpu}^{(n)}$, where $T_{\rm cpu}^{(n)}$ is the wall-clock time
    using $n$ processors.}\label{ctp-02}
\begin{tabular}{ccccccccccccccccccccccccccccc}\hline
\multirow{2}{*}{$\nu$}&\multicolumn{4}{c}{parareal}
&&\multicolumn{4}{c}{MGRiT}&&\multicolumn{4}{c}{ParaDiag-II (B-E)}
&&\multicolumn{4}{c}{ParaDiag-II (TR)} \\ \cline{2-5} \cline{7-10} \cline{12-15} \cline{17-20}
~& $16$ & $32$ & $64$ & $128$ && $16$ & $32$ & $64$ & $128$
&& $16$ & $32$ & $64$ & $128$ && $16$ & $32$ & $64$ & $128$ \\ \hline

 $10^0$    & 1.94 & 2.28 & 2.49 & 2.69 && 2.93 & 4.33 & 5.28 & 7.00 && 3.96 & 7.54 & 15.20 & 29.03 && 3.90 & 7.73 & 15.13 & 28.47 \\ 
$10^{-1}$     & 1.95 & 2.30 & 2.53 & 2.80 && 3.03 & 4.72 & 5.57 & 7.41 && 4.03 & 7.86 & 15.53 & 30.93 && 3.89 & 7.75 & 15.03 & 29.00 \\ 
 $10^{-2}$   & 1.96 & 2.32 & 2.56 & 2.83 && 2.99 & 4.58 & 5.62 & 7.35 && 3.90 & 7.76 & 15.23 & 29.00 && 3.91 & 7.77 & 15.20 & 28.84 \\ 
 $10^{-3}$   & 1.97 & 2.33 & 2.57 & 2.85 && 2.99 & 4.58 & 5.61 & 7.37 && 3.87 & 7.33 & 15.23 & 29.85 && 3.87 & 7.72 & 15.14 & 28.63 \\ 
 $10^{-4}$  & 1.95 & 2.31 & 2.55 & 2.82 && 2.95 & 4.55 & 5.61 & 7.33 && 3.89 & 7.74 & 15.17 & 29.57 && 3.91 & 7.72 & 15.29 & 29.95 \\ 
 $10^{-5}$ & 1.94 & 2.30 & 2.54 & 2.81 && 2.99 & 4.46 & 5.64 & 7.35 && 3.56 & 7.75 & 15.15 & 28.37 && 3.89 & 7.72 & 15.15 & 29.26 \\ \hline
\end{tabular}
\end{table}
Regarding the parallel efficiency measured by $T_{\rm cpu}^{(4)} /
(32\times T_{\rm cpu}^{(128)})$ \cite{ST96} ($T_{\rm cpu}^{(n)}$ is
the wall-clock time using $n$ processors), the average parallel
efficiency for ParaDiag-II (B-E) {is 92.06\%, for} ParaDiag-II
  (TR) it is 90.70\%, while it is only 22.82\% for MGRiT and 8.75\%
  for parareal.

\subsubsection{ParaDiag-II -- Parareal Variant}

 The second way to use  ParaDiag {within a stationary iteration}  is based on
formulating the coarse-grid-correction (CGC) procedure of the parareal
algorithm{\cite{LMT01,GV07}} as an all-at-once system and
{applying} ParaDiag {to it}.  The parareal algorithm is an
iterative {PinT} algorithm, {based on the updating formula}
\begin{equation}\label{eq2.15}
  U_{n+1}^{k}={\mathcal{F}^J(\Delta t, U_n^{k-1})+\mathcal{G}(\Delta T, U_n^{k})}
  -\mathcal{G}(\Delta T, U_n^{k-1}),~n=0,1,\dots, N_t-1,
\end{equation}
where $\mathcal{G}$ and $\mathcal{F}$ are called coarse and fine
propagator, specified by two time-integrators. The quantity
$\mathcal{F}^J\left(\Delta T,U^{k-1}_n\right)$ denotes a value
calculated by applying successively $J$ steps of the fine propagator
$\mathcal{F}$ to the differential equations with initial value
$U^{k-1}_n$ and the fine step size $\Delta t$. The integer
$J=\frac{\Delta T}{\Delta t}\geq2$ is called {the} {\em coarsening
  ratio}.  Let
$$
  b_{n+1}^{k-1}{:=}\mathcal{F}^J(\Delta t, U_n^{k-1})
  -\mathcal{G}(\Delta T, U_n^{k-1}).
$$
Then, the parareal algorithm is $U_{n+1}^{k}=\mathcal{G}(\Delta T,
U_n^{k})+b_{n+1}^{k-1}$. This is the so called CGC, which is a
sequential procedure and is often the bottleneck of the parallel
efficiency.  In \cite{W18}, the author proposed an idea {to
  parallelize the} CGC: suppos{ing} we have to solve an
initial-value problem
$$
\dot{U}(t)+f(U(t))=0, \quad U(0)=U_0,
$$
we apply $\mathcal{G}$ to a slightly {\em wrong} problem{, namely}
$$
  \dot{U}(t)+f(U(t))=0,~U(0)=\alpha U(T),
$$
where $\alpha\in(0, 1)$ is a free parameter. We use the linear case
$f(U)=AU$ to illustrate the details of {ParaDiag-II based on the
  parareal algorithm} (for {the} nonlinear case, see
\cite{W18}). {We also use for simplicity} Backward-Euler for
$\mathcal{G}$.  Let $\widetilde{U}_{n+1}{:=}\mathcal{F}^J(\Delta t,
U_n^{k-1})$.
The quantity $\mathcal{G}(\Delta T, U_n^{k-1})$   computed from the previous iteration is
$$
 \mathcal{G}(\Delta T, U_n^{k-1})=
 \begin{cases}
 \alpha(I_x+\Delta TA)^{-1}U_{N_t}^{k-1}, &n=0,\\
 (I_x+\Delta TA)^{-1}U_{n}^{k-1}, &n=1,2,\dots, N_t-1.
 \end{cases}
$$
{Note that all} the $N_t$ quantities
$\{\widetilde{U}_{n}\}_{n=1}^{N_t}$ and $\{\mathcal{G}(\Delta T,
U_n^{k-1})\}_{n=0}^{N_t-1}$ can be computed simultaneously {in
  parallel}. Hence,
$b_{n+1}^{k-1}=\widetilde{U}_{n+1}-\mathcal{G}(\Delta T,
U_n^{k-1})$. The parareal algorithm \eqref{eq2.15} can be rewritten {as}
\begin{equation*}
  (I_x+\Delta TA)U_{n+1}^{k}= U_n^{k}+(I_x+\Delta TA)b_{n+1}^{k-1}\quad
  \Longrightarrow \quad
  \frac{U_{n+1}^{k}- U_n^{k}}{\Delta T}+AU_{n+1}^{k}=( {\Delta T}^{-1}I_x+A)b_{n+1}^{k-1},
\end{equation*}
where $U_0^{k}=\alpha U_{N_t}^{k}$, which can be represented as
 {\small\begin{equation*}
\left(\underbrace{\frac{1}{\Delta T}\begin{bmatrix}
1 & & &-\alpha\\
-1 &1 & &\\
&\ddots &\ddots &\\
& &-1 &1
\end{bmatrix}\otimes I_x}_{=C_1^{(\alpha)}\otimes I_x}+\underbrace{\begin{bmatrix}
A & & &\\
  &A & &\\
&  &\ddots &\\
& &  &A
\end{bmatrix}}_{=I_t\otimes A}\right)\underbrace{\begin{bmatrix}
U_1^{k}\\
U_2^{k}\\
\vdots\\
U_{N_t}^{k}\\
\end{bmatrix}}_{={\bm u}^{k}}=
\underbrace{\begin{bmatrix}
({\Delta T}^{-1}I_x+A)\widetilde{U}_1- \alpha {\Delta T}^{-1}U_{N_t}^{k-1}\\
({\Delta T}^{-1}I_x+A)\widetilde{U}_2- {\Delta T}^{-1}U_{1}^{k-1}\\
\vdots\\
({\Delta T}^{-1}I_x+A)\widetilde{U}_{N_t}- {\Delta T}^{-1}U_{N_t-1}^{k-1}\\
\end{bmatrix}}_{={\bm b}^{k}}.
\end{equation*}}
This {problem is now precisely of} the form \eqref{eq2.10} {for
  $\theta=1$, i.e. the Backward-Euler method,} and the solution ${\bm
  u}^k$ can be {obtained using ParaDiag-II} (cf. \ref{eq2.12}). The
convergence rate of this {ParaDiag-II parareal variant} is
{summarized in the following theorem}.
\begin{theorem}[see \cite{W18}]\label{the1.3}
   Let $\rho_{\rm SinT-CGC}$ be the convergence factor of the parareal
   algorithm with sequential-in-time CGC (i.e., the classical parareal
   algorithm) and $\rho_{\rm PinT-CGC}$ be the convergence factor with
   parallel-in-time CGC. Then, there exists some threshold $\alpha^*$
   of the parameter $\alpha$, such that
   $$
     \rho_{\rm PinT-CGC}=\rho_{\rm SinT-CGC}, ~\text{if ~} 0<\alpha\leq \alpha^*.
   $$
   In particular, for linear {systems of} ODEs $\dot{U}(t)+AU(t)=f$ with
   $\sigma(A)\subset[0,\infty)$, i.e., all the eigenvalues of $A$ are
     non-negative real numbers, if we choose for $\mathcal{G}$ the
     {\rm Backward-Euler} method and for $\mathcal{F}$ {an} L-stable
     time-integrator (e.g., the Radau IIA methods and the Lobatto IIIC
     methods), it holds that $\alpha^*\approx0.3$.
\end{theorem}
This implies that if $\alpha$ does not exceed the threshold
$\alpha^*$, the {ParaDiag-II} parareal algorithm has the same
convergence rate as the classical parareal algorithm.

\begin{figure}
\centering
\includegraphics[width=0.49\textwidth]{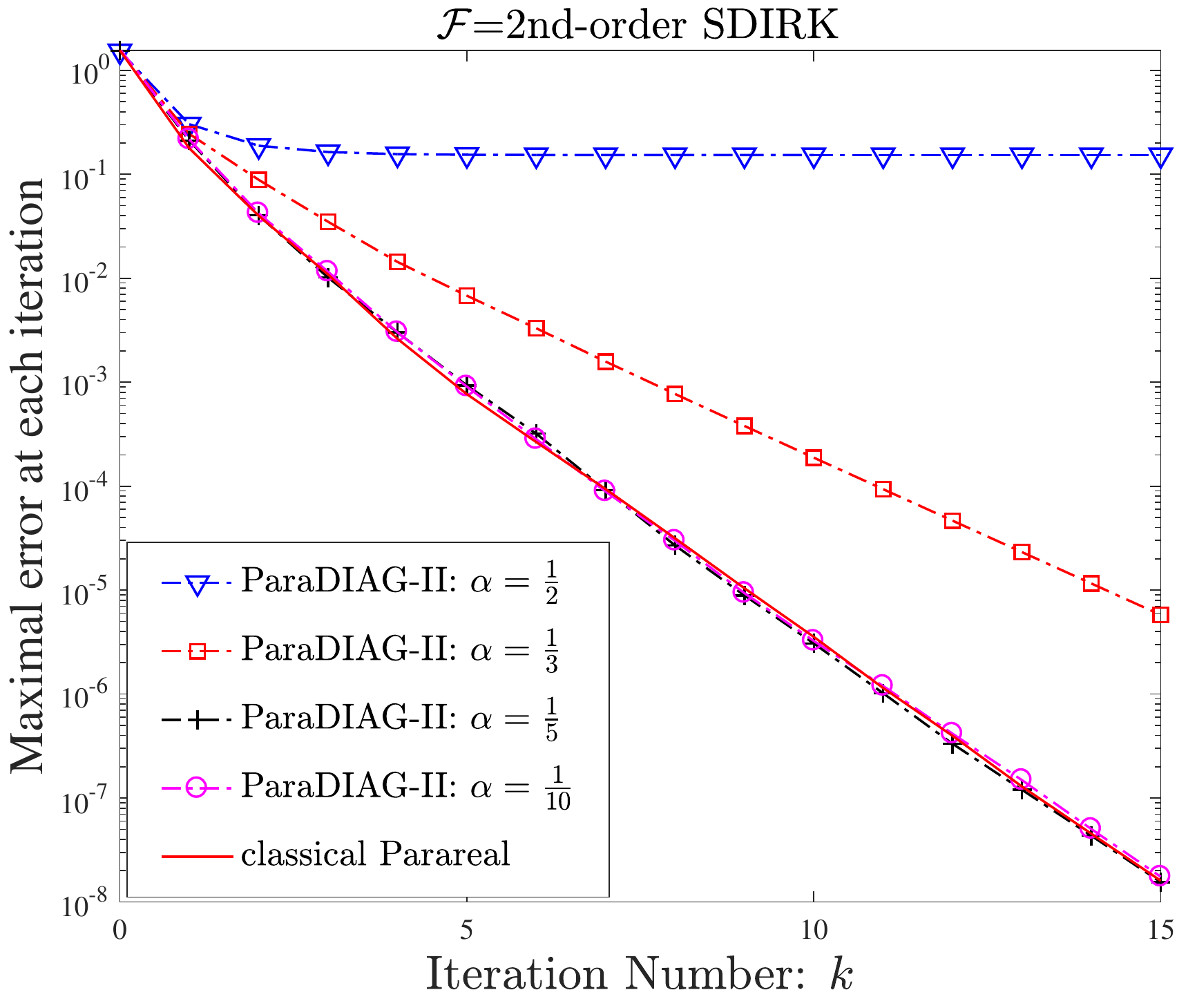}
\includegraphics[width=0.49\textwidth]{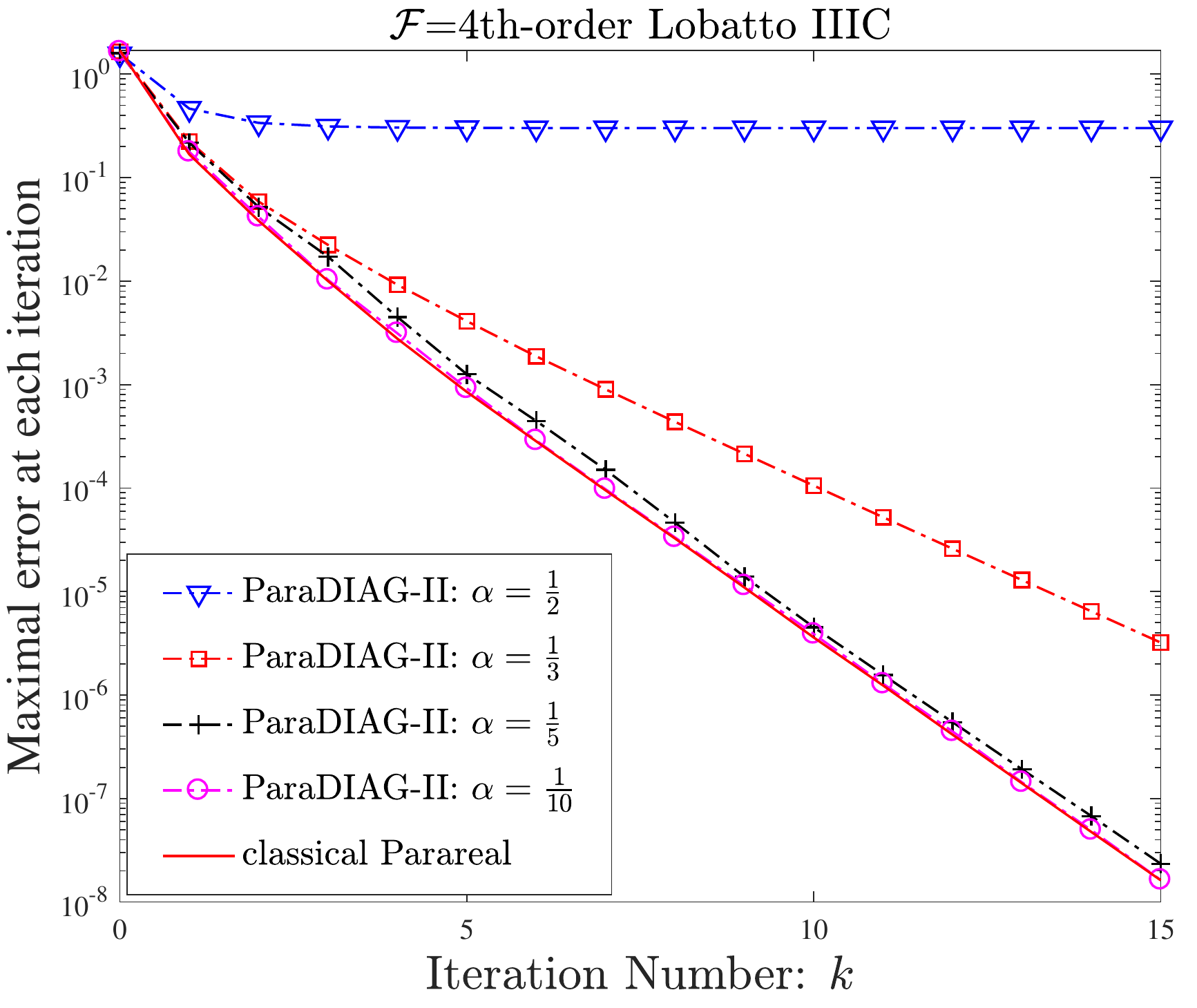}	
\caption{\em For $\nu=10^{-1}$, $T=4$, $\Delta x=\frac{1}{64}$,
  $\Delta T=\frac{1}{16}$ and $J=32$, {the measured error of
  ParaDiag-II parareal compared to classical parareal.}}
  \label{fig1.5}
\end{figure}

We {provide the} Matlab code
\texttt{ParaDiag\_V2\_Parareal\_for\_ADE} to test the convergence of
the ParaDiag-II parareal algorithm. The code includes a function
\texttt{choose\_F}, which provides 4 choices for the
$\mathcal{F}$-propagator: the Backward-Euler method, the 2nd-order
SDIRK (Singly Diagonally Implicit Runge-Kutta) method, the 3rd-order
Radau IIA method and the 4th-order Lobatto IIIC method. The interested
reader can add more choices for $\mathcal{F}$ in this function.
Moreover, we deal with a single step of the $\mathcal{F}$-propagator
by a function \texttt{Pro\_F}.  The diagonalization procedure is still
implemented via the \texttt{fft} and \texttt{ifft} commands. Starting
from a random initial guess, the error at each iteration of the new
ParaDiag-II parareal algorithm is shown in Figure \ref{fig1.5}.
The ParaDiag-II parareal algorithm converges as fast as the classical
parareal algorithm when $\alpha\leq \frac{1}{5}$.  We mention that
there is also a{n} MGRiT variant of ParaDiag-II \cite{WZ19}, which
uses {a} different {\em head-tail} coupled condition together with
the diagonalization technique for {a} parallel CGC procedure.

\section{{ParaDiag-II -- Krylov Variant (for wave equations)}}\label{sec3}

{It is a longstanding challeng{ing} task to design efficient
  PinT algorithms for wave propagation problems.  The ParaDiag-II WR
  algorithm \cite{GW19} can handle such problems, with rapid, robust
  and analyzable convergence rate. (Both the Parareal variant
  \cite{W18} and MGRiT variant \cite{WZ19} of ParaDiag-II can NOT
  handle wave equations.)  The WR variant of ParaDiag-II is an
  algorithm used within a stationary iteration.  Here, we {present}
  a ParaDiag-II {variant to be used as preconditioner for a} Krylov subspace
  method, and {which is also} applicable to wave equations. {A
    further} advantage of the new variant is that it can also
  efficiently handle optimal control problems of wave equations as
  described in Section \ref{sec4}, while {currently} some basic
  tools are lack{ing} (at least at the moment) to handle such
  optimal control problems via ParaDiag-II WR. }

{The idea below was first} introduced by McDonald, Pestana and
Wathen in \cite{MPW18} for parabolic problems, but here we {show that
  a key modification makes it a good solver also for} wave propagation
problems.  We consider the linear wave equation
\begin{equation}\label{eq3.1}
\begin{cases}
u_{tt}-\Delta u=f, &\tn{in} \Omega\times(0,T),\\
u=0,&\tn{on} \partial\Omega\times(0, T),\\
u(\cdot, 0)=u_0,~ u_t(\cdot, 0)=u_1,&\tn{in} \Omega,
\end{cases}
\end{equation}
where $\Omega\subset \mathbb{R}^d$ with $d\ge 1$ is the space domain, $u_0$ and $u_1$ are given compatible initial conditions and $f$ is a
given source term.  We discretize \eqref{eq3.1} in time by the
implicit leap-frog finite difference scheme \cite{LLX2015} (but other
schemes can be adopted as well), which was shown to be unconditionally
stable without imposing the restrictive Courant--Friedrichs--Lewy
(CFL) condition on spatial and temporal mesh sizes.  Similar to the
advection-diffusion equation, we can represent the space and time
discretizations by an all-at-once system,
\begin{subequations}
  \begin{equation}\label{eq3.2}
    {\bm A}{\bm u}:=  \left(B_1\otimes I_x +B_2\otimes A \right) {\bm u}={\bm b},
  \end{equation}
  where $A\in\mathbb{R}^{N_x\times N_x}$ is the discrete matrix of
  {the negative Laplacian} $-\Delta$, and
  \begin{equation} \label{eq3.2b}
    \begin{split}
	&B_1=\frac{1}{\Delta t^2}\begin{bmatrix}
	1 & & & &\\
	-2 &1 & & &\\
	1 &-2 &1 & &\\
	&\ddots &\ddots &\ddots &\\
	& &1 &-2 &1
	\end{bmatrix},\quad
	B_2 =\frac{1}{2}\begin{bmatrix}
	1 & & & &\\
	0 &1 & & &\\
	1 &0 &1 & &\\
	&\ddots &\ddots &\ddots &\\
	& &1 &0 &1
	\end{bmatrix}\in\mathbb{R}^{N_t\times N_t}.
    \end{split}
  \end{equation}
\end{subequations}
The idea in \cite{MPW18} for solving \eqref{eq3.2} is to construct
{a} circulant block preconditioner ${\bm P}=C_1\otimes I_x
+C_2\otimes A$, obtained by replacing the two Toeplitz matrices $B_1$
and $B_2$ in \eqref{eq3.2b} by {the} two Strang circulant matrices
\begin{equation}\label{eq3.3}
C_1=\frac{1}{\Delta t^2}\begin{bmatrix}
1 & & &1 &-2\\
-2 &1 & & &1\\
1 &-2 &1 & &\\
&\ddots &\ddots &\ddots &\\
& &1 &-2 &1
\end{bmatrix},  ~
C_2 =\frac{1}{2}\begin{bmatrix}
	1 & & &1  &\\
	0 &1 & & &1\\
	1 &0 &1 & &\\
	&\ddots &\ddots &\ddots &\\
	& &1 &0 &1
	\end{bmatrix}.
\end{equation}
Unfortunately, as {we will see later in} Table \ref{T4A}, this
preconditioner does not achieve satisfactory convergence rates for
wave equations {of the form} \eqref{eq3.1}, in contrast to parabolic
equations for which it was designed in \cite{MPW18}.

{The new idea of ParaDiag-II is to} use a generalized
preconditioner ${\bm P}_{\alpha}= \frac{1}{\Delta
  t^2}C^{(\alpha)}_1\otimes I_x + C^{(\alpha)}_2\otimes A$ by
replacing $B_1$ and $B_2$ by $\alpha$-circulant matrices (with
$\alpha\in(0,1]$ again a free parameter), where
\begin{equation}\label{eq3.4}
C^{(\alpha)}_1=\frac{1}{\Delta t^2}\begin{bmatrix}
1 & & &\alpha &-2\alpha\\
-2 &1 & & &\alpha\\
1 &-2 &1 & &\\
&\ddots &\ddots &\ddots &\\
& &1 &-2 &1
\end{bmatrix},  ~
C^{(\alpha)}_2 =\frac{1}{2}\begin{bmatrix}
	1 & & &\alpha &\\
	0 &1 & & &\alpha\\
	1 &0 &1 & &\\
	&\ddots &\ddots &\ddots &\\
	& &1 &0 &1
	\end{bmatrix}.
\end{equation}
According to Lemma \ref{lem1}, these two $\alpha$-circulant matrices
$C^{(\alpha)}_1$ and $C^{(\alpha)}_2$ can be simultaneously
diagonalized as $C^{(\alpha)}_{1,2}=V D_{1,2}V^{-1}$ and thus for an
input vector $r$ the inversion computation of ${\bm P}_\alpha^{-1}r$
can be {performed by ParaDiag} (cf. \eqref{eq1.3}): let
$D=I_x+\frac{\Delta t^2}{2}A$ and $D=Q{\rm diag}(\lambda_1,\dots,
\lambda_{N_x}) Q^\T$ be the spectral decomposition of $D$ with an
orthogonal matrix $Q$ and a real diagonal matrix ${\rm
  diag}(\lambda_1,\cdots,\lambda_{N_x})$ including all the sorted
(increasing) eigenvalues.  We have the following result for the
spectrum of the preconditioned matrix {$\CP_\alpha^{-1}{\bm A}$}.
\begin{theorem}[see \cite{WL20b}]\label{the3.1}
  The eigenvalues of the matrix {$\CP_\alpha^{-1}{\bm A}$} are explicitly
  given by
  \begin{equation*}
	\sigma({\CP_\alpha^{-1}{\bm A}})=
	\{\underbrace{1,1,\dots,1}_{(N_t-2)N_x}\}\cup  \left\{\frac{1}{1-\alpha e^{\pm \Ii N_t\theta_j}}\right\}_{j=1}^{N_x}{,}
  \end{equation*}
  where $\theta_j:=\arctan\left(\sqrt{\lambda_j^2-1}\right)\in
  (0,\pi/2)$.  Moreover, we further have the estimates:
  \begin{enumerate}
    \item If $\alpha\in (0,1)$, then
      $$
	\sigma({\CP_\alpha^{-1}{\bm A}})\subset \mathbb{A_\alpha}:=\left\{z\in\IC:\frac{\alpha}{1+\alpha} \le |z-1|\le \frac{\alpha}{1-\alpha} \right\}.
      $$
    \item If $\alpha=1$, then
      $$
	\sigma({\CP_\alpha^{-1}{\bm A}}) =
	\{1\}\cup  \left\{\frac{1}{2} \pm \frac{1}{2}\Ii\cot\left(\frac{N_t\theta_j}{2}\right) \right\}_{j=1}^{N_x}.
      $$
  \end{enumerate}
\end{theorem}
We provide {the} Matlab code \texttt{ParaDiag\_V2\_GMRES\_LinearWave\_2D} to
solve a 2D {wave equation} example with
$$
T=2, u_0(x,y)=\sin(\pi x)\sin(\pi y) ,~ u_1(x, y)=\sin(\pi x)\sin(\pi y),~ f=(1+2\pi^2)\sin(\pi x)\sin(\pi y)e^t,
$$
where the exact solution is $u(x,y,t)=\sin(\pi x)\sin(\pi y)e^t$.
Here we choose a zero initial guess and a stopping tolerance
$\texttt{tol}=10^{-10}$ based on the reduction in relative residual
norms.  The complex-shifted systems in Step-(b) are solved by MATLAB's
sparse direct solver.  We will measure the discrete
$L^\infty((0,T);L^2(\Omega))$ error norms of the numerical
approximation, and then estimate the experimental order of accuracy by
calculating the logarithmic ratio of the approximation errors between
two successively refined meshes, i.e.,
$$
\mathrm{Order} = \log_2\left(\frac{\textrm{Error} (h,\tau)}{\textrm{Error} (2h,2\tau)}\right),
$$
which should be close to 2 for second-order accuracy.  As we can see
from Table \ref{T4A}, the iteration numbers for the preconditioner
with {the original choice of $\alpha=1$} grow dramatically
{when} the mesh is refined.  {This is much better with the
  smaller choice $\alpha=0.1$ in the new ParaDiag-II algorithm, where
  we observe only 3 iterations.}  The CPU times also show the expected
quasilinear time complexity of ParaDiag-II.
\begin{table}
	\centering
	\caption{ParaDiag-II -- GMRES    for two values of the parameter $\alpha$}
	\begin{tabular}{|c|cc|cc||cc|cc|cc|cc|cc|}\hline
		&\multicolumn{4}{c|}{$\alpha=1$}&\multicolumn{4}{c|}{$\alpha=0.1$}
		\\
		\hline
		$(N_x,N_t)$& Error &Order & Iter &CPU & Error &Order & Iter &CPU    \\ \hline
		(32,32,33)&		 7.17E-03&	 1.9 &	 3&	 0.07           &7.17E-03&	 1.9 &	 3&	 0.04 \\
		(64,64,65)&		 1.86E-03&	 1.9 &	 7&	 0.57           &1.86E-03&	 1.9 &	 3&	 0.31 \\
		(128,128,129)&	4.74E-04&	 2.0 &	 37&	 24.25  & 4.74E-04&	 2.0 &	 3&	 2.17 \\
		(256,256,257)& &&$>$50&                                   & 1.20E-04&	 2.0 &	 3&	 21.02 \\ \hline
	\end{tabular}
	\label{T4A}
\end{table}

\section{{ParaDiag-II -- Krylov Variant (for optimal control of {the} wave equation)}}\label{sec4}

The Krylov variant of ParaDiag-II {can also be} used to handle optimal control problems
of the wave equation \cite{WL20b}, {by applying ParaDiag as a
  preconditioner for} the discrete saddle-point system within the
framework of Krylov subspace methods.  Let $\Omega\in \mathbb{R}^d$
with $d\ge 1$ be {a bounded and open domain with Lipschitz boundary},
and $[0, T]$ be the time window of interest with $T>0$.  We consider a
distributed optimal control problem of minimizing a tracking-type
quadratic cost functional,
 \begin{subequations}
 	\begin{equation}\label{eq2.1a}
 \min_{u,\tilde{u}}	\mathcal{L} (u ,u) :=\frac{1}{2}\|u-g\|^2_{L^2(\Omega\times(0,T))}  +\frac{\gamma}{2}\|\tilde{u}\|^2_{L^2(\Omega\times(0,T))},
 	\end{equation}
 	subject to a linear wave equation with initial- and
        boundary conditions
 	\begin{equation}\label{eq2.1b}
 	\begin{cases}
 	u_{tt}-\Delta y=f+\tilde{u}, &\tn{in} \Omega\times(0,T),\\
	 u=0, &\tn{on} \partial\Omega\times(0, T),\\
 	u(\cdot, 0)=u_0,\quad u_t(\cdot, 0)=u_1,&\tn{in} \Omega,
 	\end{cases}
 	\end{equation}
\end{subequations}
where $\tilde{u}\in L^2$ is the distributed control, $g\in L^2$ is the
desired tracking trajectory or observation data and $\gamma> 0$
is  the cost weight or regularization parameter.  {The first-order optimality system
 of \eqref{eq2.1a}-\eqref{eq2.1b} is }
 \begin{equation}\label{eq2.2}
 \begin{aligned}
 \begin{cases}
 u_{tt}-\Delta u-\frac{1}{\gamma}p=f,\quad &\tn{in} \Omega\times(0,T),\qquad y=0, \tn{on} \partial\Omega\times(0, T),\\
 u(\cdot, 0)=u_0,\quad u_t(\cdot, 0)=u_1, &\tn{in} \Omega,\\
 p_{tt}-\Delta p+u=g, &\tn{in}\Omega\times(0,T),\qquad p =0, \tn{on}\partial\Omega\times(0, T),\\
 p(\cdot, T)=0,\quad p_t(\cdot, T)=0, &\tn{in}\Omega,
 \end{cases}
 \end{aligned}
 \end{equation}
 {where we have eliminated the control variable $\tilde{u}$ from the
   optimality condition $\gamma \tilde{u}-p=0$ in \eqref{eq2.2},
   leading to a reduced optimality system regarding only $u$ and $p$.

By using the implicit leap-frog finite difference scheme
\cite{LLX2015} we get the discrete saddle-point system
\begin{equation*}
 \widehat{\bm A}
 \begin{bmatrix}{\bm u}\\ {\bm p}\end{bmatrix}:=\left(
 \begin{bmatrix}
 B_1 &-\frac{\Delta t^2\hat{I}_t}{\gamma}\\
 \Delta t^2\check{I}_t &B_1^\T
 \end{bmatrix}\otimes I_x+\frac{\Delta t^2}{2}\begin{bmatrix}
 B_2 &\\
 &B_2^\T
 \end{bmatrix}\otimes A \right)
 \begin{bmatrix}{\bm u}\\ {\bm p}\end{bmatrix}=
 \begin{bmatrix} {\bm f}\\ {\bm g}\end{bmatrix},
\end{equation*}
where $\hat{I}_t=\text{diag}(\frac{1}{2}, 1,\dots,1)$,
$\check{I}_t=\text{diag}(1,\dots,1,\frac{1}{2})\in\mathbb{R}^{N_t\times
  N_t}$, $A\in\mathbb{R}^{N_x\times N_x}$ is the discrete matrix of
{the negative Laplacian} $-\Delta$ and $B_{1,2}$ are the Toeplitz
matrices given by \eqref{eq3.2b}.  The idea in \cite{WL20b} for
applying the ParaDiag algorithm lies in three steps. First, we need to
balance the effect of the regularization parameter $\gamma$ via a
similarity transform
 \begin{equation*}
\underbrace{\left(\begin{bmatrix}\gamma^{\frac{1}{2}} I_t &\\ &I_t\end{bmatrix}\otimes I_x\right)  \widehat{\bm A} \left(\begin{bmatrix}\gamma^{-\frac{1}{2}} I_t &\\ &I_t\end{bmatrix}\otimes I_x\right)}_{:={\bm A}}
 \begin{bmatrix}\gamma^{\frac{1}{2}}{\bm u}\\ {\bm p}\end{bmatrix}= \begin{bmatrix}\gamma^{\frac{1}{2}} {\bm f}\\ {\bm g}\end{bmatrix},
 \end{equation*}
 where $
{\bm A}=\begin{bmatrix}
B_1 &-\frac{\Delta t^2\hat{I}_t}{\sqrt{\gamma}}\\
\frac{\Delta t^2\check{I}_t}{\sqrt{\gamma}} &B_1^\T
\end{bmatrix}\otimes I_x+\frac{\Delta t^2}{2}\begin{bmatrix}
B_2 &\\
 &B_2^\T
\end{bmatrix}\otimes A$. Second, based on the Toeplitz structure we
propose the following block circulant preconditioner
\begin{equation*}
{\bm P}:=
\begin{bmatrix}
C_1 &-\frac{\Delta t^2I_t}{\sqrt{\gamma}}\\
\frac{\Delta t^2I_t}{\sqrt{\gamma}} &C_1^\T
\end{bmatrix}\otimes I_x+\frac{\Delta t^2}{2}\begin{bmatrix}
C_2 &\\
&C_2^\T
\end{bmatrix}\otimes A,
\end{equation*}
where $C_1$ and $C_2$ are given by \eqref{eq3.3}. Note that the
diagonal matrices $\hat{I}_t$ and $\check{I}_t$ are replaced by the
identity matrix $I_t\in\mathbb{R}^{N_t\times N_t}$. The last step is
to rewrite ${\bm P}$ as
\begin{equation}\label{eq4.3}
 {\bm P}=
\underbrace{\left(\begin{bmatrix}
C_1C_2^{-1} &-\frac{\Delta t^2(C_2^{-1})^\T}{\sqrt{\gamma}}\\
\frac{\Delta t^2C_2^{-1}}{\sqrt{\gamma}} &C_1^\T (C_2^{-1})^\T
\end{bmatrix}\otimes I_x+\frac{\Delta t^2}{2}\begin{bmatrix}
I_t &\\
 &I_t
\end{bmatrix}\otimes A\right)}_{=:\widetilde{{\bm P}}}\left(\begin{bmatrix}C_2 &\\ &C_2^\T\end{bmatrix}\otimes I_x\right).
 \end{equation}
Now, for any input vector $r$, we can compute $s={\bm P}^{-1}r$   via
  \begin{equation*}
\tilde{s}:=
\begin{bmatrix}\tilde{s}_1\\
\tilde{s}_2\end{bmatrix}=
\widetilde{{\bm P}}^{-1}r,~s=
\begin{bmatrix}
(C_2^{-1}\otimes I_x)\tilde{s}_1\\ ((C_2^{-1})^\T\otimes I_x)\tilde{s}_2
\end{bmatrix}.
\end{equation*}
Once $\tilde{s}$ is calculated, we can compute $s$ with high
efficiency by {the} fast Fourier transform (FFT).  Hence, the major
computation is to compute $\tilde{s}=\widetilde{{\bm P}}^{-1}r$.  We
now derive a special diagonalization of the matrix $\widetilde{{\bm
    P}}$ in \eqref{eq4.3}.  The reason why we split ${\bm P}^{-1}$
into two steps is that we do not have a spectral decomposition of
${\bm P}$ with {a} closed formula. For $\widetilde{{\bm P}}$, we have the
following spectral decomposition.
\begin{theorem}[{see} \cite{WL20}]\label{the4}
Let $D_1$ and $D_2$ be  the diagonal matrices consisting of the circulant matrices $C_1$ and $C_2$ and  $\mathbb{F}\in\mathbb{C}^{N_t\times N_t}$ be the discrete Fourier matrix.  The matrix  $\widetilde{{\bm P}}$ in \eqref{eq4.3} can be factorized as
\begin{subequations}
\begin{equation}\label{eq4.4a}
\begin{split}
&\widetilde{{\bm P}}=(V\otimes I_x)\left(\begin{bmatrix}\Sigma_1 &\\
&\Sigma_2\end{bmatrix}\otimes I_x+
\frac{\Delta t^2}{2}
\begin{bmatrix}
I_t &\\
&I_t
\end{bmatrix}\otimes  A\right)(V^{-1}\otimes I_x).
\end{split}
\end{equation}
where
\begin{equation}\label{eq4.4b}
\begin{split}
&V=\begin{bmatrix}
	\mathbb{F}^* &\\
	&\mathbb{F}^*
	\end{bmatrix} \begin{bmatrix}
I_t &-{\rm i} \sqrt{D_2^*D_2^{-1}}\\
{\rm i} \sqrt{D_2^*D_2^{-1}} &I_t
\end{bmatrix}, ~~V^{-1}=\frac{1}{2}V^*, \\
&\Sigma_1=
D_1D_2^{-1}+ \Ii \frac{\Delta t^2}{\sqrt{\gamma}}|D_2^{-1}|, ~~\Sigma_2=
D_1D_2^{-1}- \Ii \frac{\Delta t^2}{\sqrt{\gamma}}|D_2^{-1}|.
\end{split}
\end{equation}
\end{subequations}
\end{theorem}
Let $\widetilde{D}$ be an invertible diagonal matrix. Then, it is
clear that the factorization \eqref{eq4.4a} still holds if we replace
$V$ by $V\widetilde{D}$.  Hence the eigenvector matrix for the block
diagonalization of $\widetilde{{\bm P}}$ is not unique. A nice
property of the factorization given by \eqref{eq4.4a}-\eqref{eq4.4b}
is that the matrix $V$ is {\em optimal} in the sense that
Cond$_2(V)=1$.  According to \eqref{eq4.4b}, for any input vector
${\bm r}$ we can compute $\widetilde{{\bm P}}^{-1}{\bm r}$ by the
diagonalization technique described {in} \eqref{eq1.3}.  It was
shown in \cite{WL20} that the eigenvalues of the non-symmetric
preconditioned matrix ${\bm P}^{-1}{\bm A}$ {are} highly clustered
(the similarity transform from $\widehat{\bm A}$ to ${\bm A}$ is
important for this).

We provide {the} Matlab code \texttt{ParaDiag\_V2\_GMRES\_LinearWaveOPT\_2D} for
the {2D wave equation optimal control problem} posed on
$\Omega\times(0,T)=(0,1)^2\times(0,2)$, with the data
\begin{equation*}
\begin{split}
&u_0(x,y)=\sin(\pi x)\sin(\pi y),\quad u_1(x,y)=\sin(\pi x)\sin(\pi y),\\
&f(x,y,t)=(1+2\pi^2)e^t \sin(\pi x)\sin(\pi y)-\frac{1}{\gamma}(t-T)^2\sin(\pi x)\sin(\pi y),\\
&g(x,y,t)=(e^t+2+2\pi^2(t-T)^2)\sin(\pi x)\sin(\pi y).
\end{split}
\end{equation*}
The exact solution of the optimal control problem is
\[
u(x,y,t)=e^t\sin(\pi x)\sin(\pi y)\quad\mbox{and}\quad p(x,y,t)=(t-T)^2\sin(\pi x)\sin(\pi y).
\]
As shown in Table \ref{T2A},
\begin{table}
 	\centering
 	\caption{{Number of GMRES iterations and CPU times using the
            {ParaDiag-II} preconditioner ${\bm P}$.}}
 	\begin{tabular}{|c||cc|cc|cc|cc|cc|cc|cc|}\hline
 		$\texttt{tol}=10^{-7}$ &\multicolumn{2}{c}{$\gamma=10^{-2}$}&\multicolumn{2}{|c|}{$\gamma=10^{-4}$}
 		&\multicolumn{2}{c}{{$\gamma=10^{-6}$}}&\multicolumn{2}{|c|}{$\gamma=10^{-8}$}
 		&\multicolumn{2}{c|}{{$\gamma=10^{-10}$}}
 		\\
 		\hline
 		$(N_x,N_x,N_t)$&It&CPU&It&CPU&It&CPU&It&CPU&It&CPU    \\ \hline
 	(16,16,17)&	            5&	 0.0 &          5&	 0.0    &      5&	 0.0    &     4&	 0.0     &        4&	 0.0 \\
 	(32,32,33)&	            5&	 0.1 &          5&	 0.2    &      5&	 0.1    &     5&	 0.1     &        4&	 0.1 \\
 	(64,64,65)&	            5&	 0.7 &          5&	 1.1    &      5&	 0.8    &     5&	 0.8     &        4&	 0.6 \\
 	(128,128,129)&          11&	 13.9 &        5&	 7.2    &      5&	 6.7    &     5&	 6.4     &        5&	 6.7 \\
 	(256,256,257)&          17&	 226.6&        5&	 59.6   &      5&	 60.3   &     5&	 61.0    &        5&	 60.7 \\
 		\hline
 	\end{tabular}
 	\label{T2A}
\end{table}
GMRES {preconditioned with the
  ParaDiag-II} preconditioner ${\bm P}$ converges very fast and {is
robust with} respect to the possibl{y} very small regularization parameter
$\gamma$.

\section{ParaDiag-II: A General Theoretical Result (\red{\bf New Progress})}\label{sec5}

{For the all-at-once system \eqref{eq1.1}, {using ParaDiag as
    a stationary iterative solver corresponds to the iteration}
\begin{equation}\label{eq5.1}
{\bm P}_\alpha\Delta {\bm u}^k={\bm r}^k, ~{\bm u}^{k+1}={\bm u}^k+\Delta{\bm u}^k,~{\bm r}^k:={\bm b}-{\bm A}{\bm u}^k, 
\end{equation}
where ${\bm P}_\alpha$ is the block $\alpha$-circulant matrix defined
by \eqref{eq1.2}. In \cite{MPW18,LN20}, such a ${\bm P}_\alpha$ (with
$\alpha=1$ in \cite{MPW18}) was used as a preconditioner for 
Krylov subspace solvers as we explained in Section \ref{sec3}.  But in
\cite{WL20b}, we show that the stationary iteration itself performs
very well for {both parabolic} and hyperbolic problems.  In
particular, for the implicit leapfrog scheme it was proved {that}
\begin{equation}\label{eq5.2} 
\rho({\bm I}-{\bm P}_\alpha^{-1}{\bm A})\leq\frac{\alpha}{1-\alpha},~\alpha\in(0, 1), 
\end{equation}
where the upper bound only depends on $\alpha$. Actually, this is the
same statement {as in} Theorem \ref{the3.1}.  For the Krylov variants of
ParaDiag-II introduced in Section \ref{sec3}, the eigenvalue
distribution of the preconditioned matrix ${\bm P}^{-1}_\alpha{\bm A}$
  is also an important {issue}, even though a
clustering of the eigenvalues does not necessarily imply fast
convergence of the algorithm. Considerable efforts have been devoted
to exploring the spectrum of ${\bm P}^{-1}_\alpha{\bm A}$ (or ${\bm
  I}-{\bm P}^{-1}_\alpha{\bm A}$), and this leads to many case-by-case
studies depending on the {time-integrator used}. To name a few, we
mention the work in \cite{LN20} for the implicit Euler method,
\cite{WL20b} for the implicit leap-frog method, \cite{WZ21} for the
two-stage singly diagonal implicit RK method and \cite{WZ21b} for the
BDF method with order up to 6.  The analysis in these {references} is very
technical and heavily depends on the special property of the
time-integrator, e.g., sparsity, Toeplitz structure and diagonal
dominance of the time-discretization matrix.

It is {therefore justified to ask the question:} ``{\em In general,
  under what conditions {does} the iterative algorithm
  \eqref{eq5.1} converge rapidly and robustly?}''  We recently proved
the following unified results \cite{WZZ21}:
\begin{theorem}
  For {an} initial-value problem $U'+AU=F$ with
  $A\in\mathbb{C}^{m\times m}$ and spectrum
  $\sigma(A)\subset\mathbb{C}^+$, suppose $A$ is diagonalizable as
  $A=PD_AP^{-1}$. For any one-step time-integrator
\begin{equation}\label{eq5.3}
U_{n+1}+\mathcal{R}(\Delta tA)U_{n}=\tilde{F}_n, ~n=0, 1,\dots, N_t-1,  
\end{equation}
the error at the $k$-th iteration \eqref{eq5.1}, denoted by ${\bm{err}}^{k}={\bm u}^{k}-{\bm u}$,
satisfies 
\begin{equation}\label{eq5.4}
\|(I_t\otimes P){\bm{err}}^{k+1}\|_{\infty}\leq\frac{\alpha}{1-\alpha}\|(I_t\otimes P){\bm{err}}^{k}\|_{\infty},~\forall k\geq1,
\end{equation} 
provided the time-integrator is stable in the sense 
 \begin{equation}\label{eq5.5}
 |\mathcal{R}(\Delta t \lambda)|\leq1,  {\forall \lambda\in\sigma(A)}{.}
 \end{equation}
Therefore, for one-step time-integrators the iteration \eqref{eq5.1}
converges linearly if $\alpha\in(0,\frac{1}{2})$.
\end{theorem}
For the one-step time-integrator \eqref{eq5.3}, $\mathcal{R}(\Delta
A)$ is the increment matrix deduced from the stability function. For
example, for a general implicit {$s$ stage} RK method specified by
the Butcher tableau
$$
 \begin{array}{r|l}
c &{{\Theta}}    \\
\hline
&b^\top
\end{array},
$$
{the increment matrix is given by}
  $$
  \mathcal{R}(\Delta tA)= I_x-b^\top\otimes (\Delta tA)(I_s\otimes I_x+\Theta\otimes(\Delta tA))^{-1}({\bm 1}\otimes I_x),
  $$
   where $I_s\in\mathbb{R}^{s\times s}$ is an identity matrix and ${\bm 1}=(1,1,\dots, 1)^\top\in\mathbb{R}^s$.
  
\begin{theorem}\label{the5.2}
For {an} initial-value problem 
$U'+AU=F$ with  $A\in\mathbb{C}^{m\times m}$ and  spectrum  $\sigma(A)\subset\mathbb{C}^+$, suppose  $A$ is diagonalizable as $A=PD_AP^{-1}$. For any linear multistep  method \begin{equation}\label{eq5.6}
 {\sum}_{j=0}^r a_j U_{n+r-j} + \Delta t {\sum}_{j=0}^r b_j A  U_{n+r-j}=\tilde{F}_n, \quad n=0,1,\ldots, N_t-r, 
\end{equation}
the error at the $k$-th iteration \eqref{eq5.1}, denoted by ${\bm{err}}^{k}={\bm u}^{k}-{\bm u}$,
satisfies 
\begin{equation}\label{eq5.7}
\|(I_t\otimes P){\bm{err}}^{k+1}\|_{\infty}\leq\frac{c\alpha}{1-c\alpha}\|(I_t\otimes P){\bm{err}}^{k}\|_{\infty},~\forall k\geq1,
\end{equation} 
provided the method is stable in the sense 
\begin{equation}\label{eq5.8}
p(s; z)=0~~\Longrightarrow~~
\begin{cases}
 {\rm either}~~ |s| < 1,\\
\text{\rm or}~~|s|=1~{\rm and ~it~ is~ a~ root~ of ~multiplicity~ 1},
\end{cases}
 \end{equation}
where $c\geq1$ is a constant only depending on the stability of the time-integrator,   $z=\Delta t\lambda$ is an arbitrary eigenvalue of $\Delta tA$ and  $p(s,z)$ is  the characteristic polynomials  of the $r$-step method{,}
\begin{equation}\label{eqn:char-poly}
 p(s;z) = {\sum}_{j=0}^r a_j s^{r-j} + z b_j s^{r-j}.
\end{equation}
Therefore, for multistep time-integrators the iteration \eqref{eq5.1} converges linearly if $\alpha\in(0,\frac{1}{c})$.
\end{theorem}
}

\end{document}